\documentclass[11pt]{article}

\RequirePackage{fix-cm}
\usepackage{mathptmx}
\DeclareMathAlphabet{\mathpzc}{OT1}{pzc}{m}{it}

\usepackage{bm}
\usepackage{latexsym}
\usepackage{amssymb}
\usepackage{amsfonts}
\usepackage{amsmath}
\usepackage{amsthm}
\usepackage{mathrsfs}
\usepackage{dsfont}    
\usepackage{bbold}     
\usepackage[english]{babel}
\usepackage{epsfig}
\usepackage{float}
\usepackage{subfigure}
\usepackage{psfrag}
\usepackage{graphicx}
\usepackage{epsfig}
\usepackage[margin=1.2cm]{caption}

\newtheorem{rmk}{Remark} 

\newcommand{\e}{\varepsilon}
\newcommand{\al}{\alpha}
\newcommand{\om}{\omega}

\newcommand{\oo}{\boldsymbol{\omega}} 

\newcommand{\ZZZ}{\mathds{Z}}

\newcommand{\QQQ}{\mathds{Q}}
\newcommand{\RRR}{\mathds{R}}
\newcommand{\TTT}{\mathds{T}}

\newcommand{\KK}{{\mathcal K}}
\newcommand{\II}{{\mathcal I}}
\newcommand{\QQ}{{\mathcal Q}}

\newcommand{\TT}{{\mathpzc T}}
\newcommand{\UU}{{\mathcal U}}
\newcommand{\st}{\scriptstyle}
\newcommand{\io}{\infty}
\newcommand{\Th}{\bar\Theta}

\begin{document}
\title{\textbf{Periodic and quasi-periodic attractors for the spin-orbit evolution of Mercury with a realistic tidal torque}}

\author
{\bf Michele Bartuccelli$^1$, Jonathan Deane$^1$, Guido Gentile$^2$
\vspace{2mm}
\\ \small 
$^1$ Department of Mathematics, University of Surrey, Guildford, GU2 7XH,
UK 
\\ \small
$^2$ Dipartimento di Matematica, Universit\`a Roma Tre, Roma, I-00146,
Italy
\\ \small 
E-mail: m.bartuccelli@surrey.ac.uk, j.deane@surrey.ac.uk,
gentile@mat.uniroma3.it
}
\date{}
\maketitle

\begin{abstract}
Mercury is entrapped in a 3:2 resonance: it rotates on its axis three times
for every two revolutions it makes around the Sun. It is generally accepted that this is due to the large value 
of the eccentricity of its orbit. However, the mathematical model originally introduced to study its spin-orbit evolution
proved not to be entirely convincing, because of the expression commonly used for the tidal torque.
Only recently, in a series of papers mainly by Efroimsky and Makarov, 
a different model for the tidal torque has been proposed, which has the advantages of being more realistic,
and of providing a higher probability of capture in the 3:2 resonance with respect to the
previous models. On the other hand, a drawback of the model is that the function describing
the tidal torque is not smooth and consists of a superposition of kinks, so that both analytical
and numerical computations turn out to be rather delicate: indeed, standard perturbation theory
based on power series expansion cannot be applied and the implementation of a fast algorithm
to integrate the equations of motion numerically requires a high degree of care. In this paper, 
we make a detailed study of the spin-orbit dynamics of Mercury, as predicted by the realistic model:
In particular, we present numerical and analytical results about the nature of the librations of Mercury's spin
in the 3:2 resonance. The results provide evidence that the
librations are quasi-periodic in time.
\end{abstract}

\section{Introduction} \label{sec:1}
\setcounter{equation}{0}

The reason that Mercury is entrapped in the 3:2 resonance has been extensively investigated
in the literature. It is commonly accepted that this is due to the high value of Mercury's eccentricity ($0.2056$);
however, there is no universal consensus about the mechanism by which the entrapment has occurred.

The Mercury-Sun system is usually studied as a satellite-planet system, with the satellite described as
an ellipsoidal body orbiting around its primary in a Keplerian orbit.
If $\theta$ denotes the sidereal angle, that is, the angle that the longest axis of the satellite
forms with respect to the line of apsides of the orbit, the time evolution of $\theta$ is described
by the second order ordinary differential equation
\begin{equation} \label{eq:1.1}
C \ddot\theta = \TT_z^{\rm (TRI)} + \TT_z^{\rm (TIDE)}  ,
\end{equation}
with $C$ being the maximal moment of inertia of the satellite, and where $\TT_z^{\rm (TRI)}$ and $\TT_z^{\rm (TIDE)}$
are traditionally called the triaxiality-caused (or simply triaxial) torque and the tidal torque.

There is general agreement in the literature as to the expression for $\TT_z^{\rm (TRI)}$
(see for instance Danby \cite{D}): the triaxial torque is written as an infinite
Fourier series, which is usually truncated, since only a few resonances are really relevant.
By contrast, the expression for $\TT_z^{\rm (TIDE)}$ is a much more delicate issue.
In the first paper devoted to the problem, by Goldreich and Peale \cite{GP},
the MacDonald model was mainly used for the tidal torque. This entailed a simple form for $\TT_z^{\rm (TIDE)}$,
well suited to both analytical and numerical computations. However, the results obtained with such a model
were rather disappointing: the probability of capture in the 3:2 resonance
was found to be very small (7\%), and only assuming a chaotic evolution of the
eccentricity of Mercury --- as demonstrated much later by Correia and Laskar \cite{CoLa} ---
can it become much higher (55\%). Here, the probability of capture is defined as the probability
of the satellite being trapped in a resonance when crossing it;  if one takes into account the strong variations which
the eccentricity underwent in the past, the probability is highly enhanced because multiple crossings become possible.

Moreover, tidal models such as MacDonald's,
based on a constant time lag (CTL), lead to the existence of a stable quasi-synchronous solution,
which in the case of Mercury is characterised by a spin rate $\dot\theta \approx 1.26 n$,
where $n$ is the mean
motion of Mercury; such a solution turns out to attract most of trajectories in the case of constant eccentricity.
More generally, the CTL model produces rather nonphysical results when applied to satellite-primary systems,
as the most common resonance for satellites is the synchronous one (1:1).

Recently, the physical validity of tidal models based on constant time lag was strongly questioned by Efroimsky and
Makarov \cite{E1,E2,EM,ME}. A more realistic model has been introduced by Efroimsky \cite{E2},
based on the Darwin-Kaula expansion of the tidal torque \cite{K}, which takes into account both the rheology and
the self-gravitation of Mercury. By relying on such a model, Makarov showed that the probability of capture
in the 3:2 resonance is 100\%, that is to say when Mercury crosses the 3:2 resonance it is inevitably entrapped in it \cite{M}.
Later on, Noyelles \textit{et al.}\ studied the case of non-constant eccentricity and found not only
that trapping in the 3:2 resonance is the most probable outcome of the time evolution,
but also that the trapping time is much smaller than that predicted by previous theories \cite{NFME}.
We refer to \cite{NFME} for more details and for a very clear discussion of the existing results
in the literature; in the following we shall refer to the system with the tidal model used in \cite{NFME} as the NFME model.
We note that, in the NFME model, the tidal torque is not a smooth function (it is only $C^1$) and
it appears as a superposition of kinks: this makes both the numerical and the analytical investigations rather subtle.

In general, because of the small value of the dissipation, integrating the equations of motion requires very long times.
Thus, it may be convenient for practical purposes to make some assumptions on the initial data of the system:
usually one fixes the initial velocity and considers a large sample of initial phases (say, 1000).
Therefore, a very high probability of capture (even 100\%) in a given resonance does not necessarily
imply that every trajectory ends up in that resonance, because, for that to happen, one needs
the trajectory not to have been trapped earlier by other resonances that it has crossed during its time evolution.
Since the initial condition is not known, it may be important to investigate a larger sample
of initial conditions by varying $\dot\theta$ as well as
$\theta$, randomly distributed in the phase space.
A fast numerical integration method was proposed by Bartuccelli \textit{et al.}~\cite{BDG1,BDG2},
which allows consideration of a larger number of initial data (say, 50~000). Then one can evaluate
the probability of capture in a given resonance as the fraction of trajectories which are eventually
attracted into that resonance: for the NFME model it was found that the 3:2 resonance is still 
the most probable final state, since it attracts about 42\% of the trajectories with initial conditions
$(\theta,\dot\theta)$ inside the set $[0,2\pi]\times[0,5n]$, and even more if one takes into account
only initial conditions above the 3:2 resonance (see Section \ref{sec:8}).

In this paper, we study more closely, both numerically and analytically, the nature of the attractors
for the NFME model --- something which is still missing in the literature.
Indeed, as noted in~\cite{MFD}, the tidal torque in the NFME model
``leaves little room for analytical applications" and usually very strong simplifying assumptions
on the equation are made in order to obtain analytical expressions for approximate solutions.
In fact, when one speaks of a resonance $p$:$q$, usually one simply means that the solution $\theta(t)$
to \eqref{eq:1.1} is such that $\dot\theta(t) \approx p n/q $, with $n$ being the mean motion, but
it is not obvious at all whether the solution is periodic, i.e. has frequency commensurate with $n$.

\begin{figure}[!ht]
\centering 
\vspace{-0.2cm}
\subfigure{\includegraphics[width=6.2in]{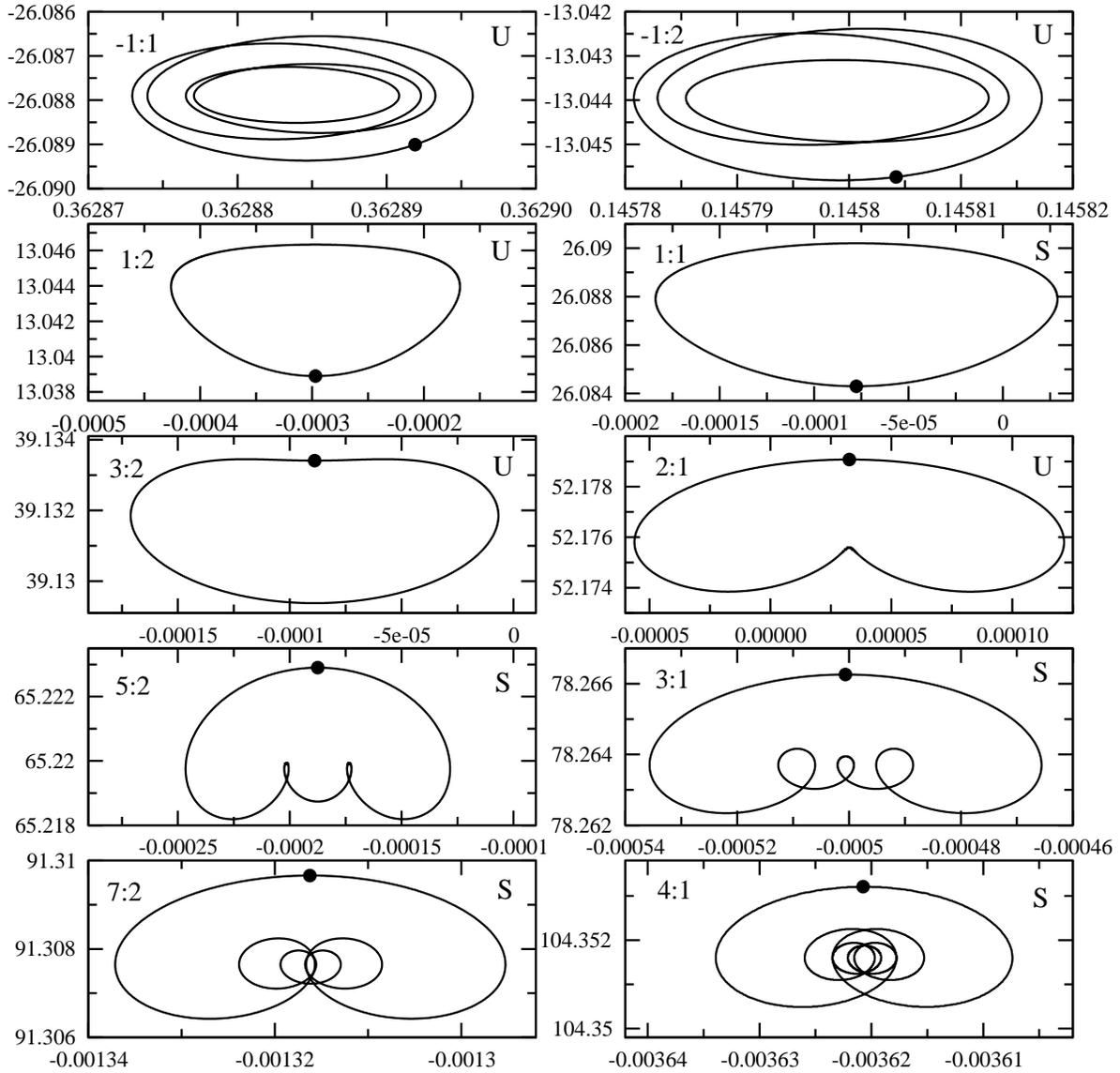}}
\caption{Periodic solutions found numerically, corresponding to the resonances $p$:$q$,
with $q=2$ and $p=-2,\ldots,8$ ($p\neq0$). The velocity $\dot\theta(t)$ is plotted versus $\theta(t)-nt$;
in each figure the dot represents the initial condition. The letters S/U denote whether a solution is stable/unstable
(as discussed in Section~\ref{sec:6}).
}
\label{figure:1}
\end{figure} 

Numerically, one finds a finite number of periodic attractors. The exact number depends
on the truncation of the triaxial torque $\TT_z^{\rm (TRI)}$.
It is important to note that, if we enlarged the number of harmonics included in  $\TT_z^{\rm (TRI)}$
with respect to the truncation used in \cite{NFME}, it is true that new periodic attractors would appear,
but in fact they would attract only a few trajectories: as a consequence,
the general scenario would remain essentially unchanged.
Using the same truncation as in \cite{NFME}, the periodic solutions in Figure \ref{figure:1}
are found numerically. Not all such solutions are stable (see Section \ref{sec:6} for details);
the unstable ones do not correspond to attractors.
In particular, the periodic solutions which are unstable include the 3:2 resonance.
This may seem a bit surprising, since the 3:2 resonance is expected to be the dominant one.

\begin{figure}[H]
\centering 
\null\hspace{0.9cm}\includegraphics[width=3.1in]{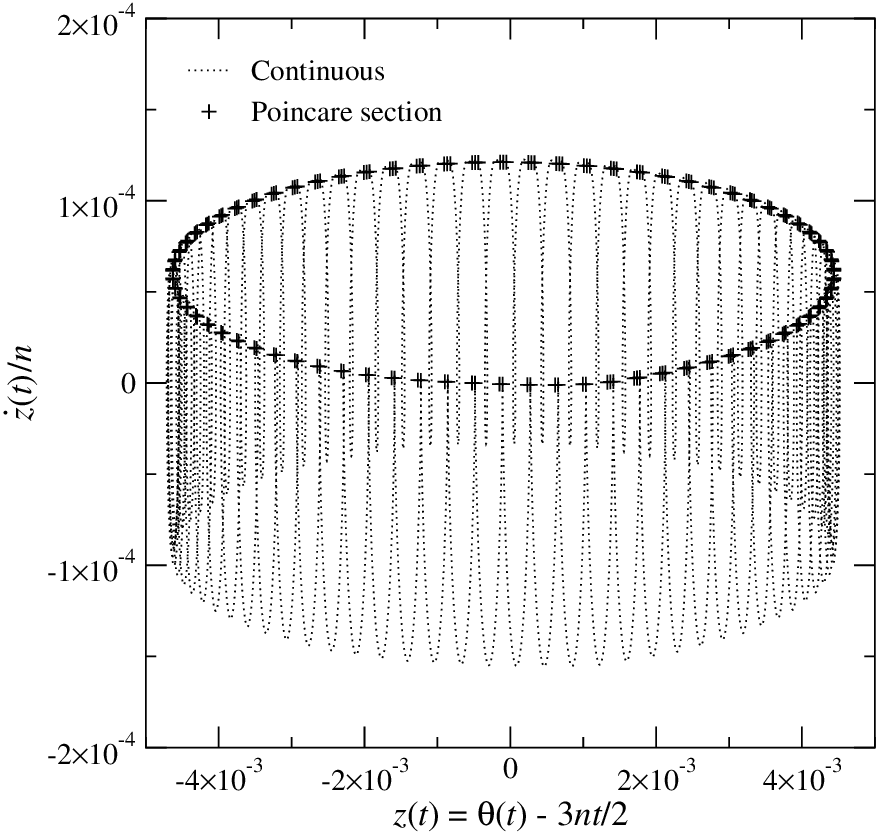}
\caption{
Quasi-periodic attractor close to the resonance 3:2, found numerically. The crosses correspond
to the Poincar\'e section, $(z(kT_0), \dot z(k T_0)/n)$, $k = 1,\ldots 200$, and the dotted line shows
the whole solution for $t = 0$ to $75 T_0$, where $T_0:=2\pi/n$.}
\label{figure:2}
\end{figure}

In fact, an attracting solution with $\dot\theta(t) \approx 3n/2$ is found numerically. However,
if we look more carefully at such a solution, we realise that it does not appear
to be a periodic solution.
More precisely, if we write $\theta(t)=3nt/2+z(t)$ and plot $\dot z(t)/n$ versus $z(t)$,
we obtain the curve in Figure \ref{figure:2} (the function $z(t)$ describes the librations of the spin rate).
So, the attracting solution has a much more complicated structure
with respect to the periodic solution with $\dot\theta(t) \approx 3n/2$ depicted in Figure~\ref{figure:1}.
Apparently, the solution is characterised by two frequencies: there is a fast oscillating motion superimposed
on a slow oscillation. This is confirmed by a Fourier Transform analysis:
the dynamics involves two frequencies $n$ and $\omega$, with $n \approx 73.9 \omega$ (see Section \ref{sec:5}).
We term such a solution a quasi-periodic attractor (see Remark \ref{rmk:6} in Section \ref{sec:5}, though).

So, from a numerical point of view, we find that the main attractor of the system \eqref{eq:1.1} does not really correspond
to what one usually means by a resonance, that is, a periodic solution with frequency commensurate with the forcing frequency.
As a matter of fact, while the mechanism by which periodic attractors appear in a periodically perturbed system is rather clear,
as follows from Melnikov's theory~\cite{GH},
the appearance of quasi-periodic solutions is much less standard and deserves further investigation.

We aim to provide an analytical description of the attractors represented in both Figures \ref{figure:1} and \ref{figure:2}.
Essentially, we shall use perturbation theory,
but with some caveats, since the tidal torque is not a smooth function and has very rapid variations. We shall see
that a few steps of perturbation theory are sufficient to provide an analytical expression for the solutions which
is in very good agreement with the numerical results. 
However, before entering into the mathematical details, let us discuss briefly --- and informally ---
what kind of solutions may be expected.

We can rewrite \eqref{eq:1.1} as
\begin{equation} \label{eq:1.2}
\ddot  \theta = -\e \, G( \theta,t) - \e \, \gamma \, F(\dot \theta) ,
\end{equation}
where $\e$ and $\gamma$ are positive numbers; the explicit expression of the functions $F$ and $G$
will be given in Section \ref{sec:2}. The number $\e$ is small, so it plays the role of a perturbation parameter.
Let us consider first the case in which $\e\neq0$ and $\gamma=0$ (conservative limit). In that case,
for $\e$ small enough, one has a quasi-integrable system and hence KAM theorem applies \cite{AKN}:
most invariant tori persist, while the resonant ones are destroyed.
Therefore, most of the solutions are quasi-periodic. This does not mean that periodic solutions do not exist:
what happens is that of each entire resonant torus only a finite number of trajectories survive in the presence
of the perturbation. When the dissipative term is also present ($\gamma \neq 0$), the scenario changes drastically:
while the periodic solutions persist, the quasi-periodic solutions disappear almost completely.
Moreover, the periodic solutions assume a pivotal role, since they become attractors. In practice,
both $\e$ and $\gamma$ are different from zero, so that the scenario may be somewhat different.  In particular,
according to the exact form of the dissipative term, quasi-periodic attractors are still possible for the full system.
For instance, this is what happens when the MacDonald torque is considered~\cite{CC,BDG0}:
the quasi-periodic solution corresponds to the quasi-synchronous solution which is also found numerically.
While no other quasi-periodic solutions are observed, analytically the existence of the quasi-synchronous solution
is rather tricky to prove. Indeed, in order to implement a KAM-like scheme, one needs to assume that the frequencies
of the solution are strongly non-resonant --- in practice Diophantine (see for instance~\cite{CC,BDG1,MNT}). 
However, the frequencies depend continuously on the parameters and the non-resonance condition is not necessarily
satisfied when the parameters are varied.

In the case of the realistic tidal model considered by Efroimsky and Makarov, stable quasi-synchronous solutions
are not possible \cite{ME}. The attracting solutions are not necessarily periodic though. Indeed, for $\e$ small enough,
periodic solutions corresponding to resonances are expected to exist and be attractive, as follows
from Melnikov's theory, but, when increasing the value of $\e$, bifurcation phenomena may occur:
a periodic solution may become unstable and a new kind of solution may appear nearby.
For the corresponding Poincar\'e map  --- or a suitable iteration of it --- the periodic solution gives a fixed point;
then the new solution looks like a curve, which emerges from the fixed point by Hopf bifurcation when
$\e$ crosses a threshold value~\cite{MC}.

A natural question is whether one can account analytically for the attracting solutions that we have just
described. As far as the solution is periodic,
one can apply Melnikov's theory for subharmonic solutions; see for instance~\cite{GBD} and references therein.
Care is in order in that case, because the function $F$ in~\eqref{eq:1.2} is only $C^{1}$
and has rapid variations, so one cannot apply perturbation theory as usually implemented
in the case of smooth functions.

A mathematical description of the solution appearing by Hopf bifurcation is more demanding.
To understand how to proceed, let us make a step back and consider the conservative system with $\gamma=0$.
In that case, as we have said, the resonant tori are destroyed. In addition, the closest tori also undergo the same fate,
so that narrow gaps appear in phase space where the resonant tori have disappeared. It is in these gaps that
the periodic orbits are located. Nevertheless, not all the motions inside the gaps are chaotic.
In fact, there are a lot of Lagrangian tori~\cite{AKN,MNT,BC}: if one considers the Poincar\'e map, the
tori look like closed curves encircling the fixed point which corresponds to the periodic solution.
Of course, such tori are not KAM tori (that is, they are not deformations of the tori of the unperturbed system),
but they still correspond to quasi-periodic solutions.

Now, let us consider what happens when the dissipation is taken into account.  If the dissipation is large
enough, then all the tori disappear and the fixed point turns out to be asymptotically stable:
such a scenario corresponds to the periodic solutions found by Melnikov's theory.
On the other hand, if the dissipation is small, the force may prevail and push away from the fixed point
(which still exists, but becomes unstable, this being an effect of the force dominating the dissipation)
and one of the curves nearby survives and is stable. 
Such a curve is the closed curve which arises by Hopf bifurcation from the fixed point.

To summarise, we find an attracting solution close to each resonance: whether such a solution is periodic
or has a more complicated structure depends on the values of the parameters. For fixed values of the
parameters, it may happen that some periodic solutions are stable (and hence are attractors)
and others are unstable (and so a quasi-periodic solution appears nearby);
in the case of Mercury-Sun, for the physical values of the parameters --- see Section~\ref{sec:2} --- we
find that the attracting solution corresponding to the 3:2 resonance is quasi-periodic. 
Both kinds of solutions can be found by perturbation theory. The periodic solutions to the full system \eqref{eq:1.2}
are obtained by starting from the periodic solutions of the unperturbed system ($\e=0$).
By contrast, the quasi-periodic solutions are not perturbations of solutions
of the unperturbed system: one has to perform first a change of variables (a `normal form', in the language of KAM theory),
which allows us to write the system as a perturbation of a new unperturbed system which is essentially a pendulum. Then
the unperturbed solution to be continued is an oscillatory solution of the pendulum.

In the forthcoming sections, we use perturbation theory to compute approximations to the periodic and quasi-periodic
solutions described above. To overcome the lack of smoothness of the vector field, we use
an iteration method based on the same idea as the Picard approximants. In principle, one should prove that the
iteration scheme converges, but we do not address this issue here. We confine ourselves to computing a few
steps of the iteration and show that the results are in very good agreement with the numerical simulations.
The problem of convergence is certainly non-trivial, particularly in the case of quasi-periodic solutions.
In fact, as pointed out previously, in KAM-like problems one usually assumes strong resonance conditions on the frequencies,
whereas, in the case of dissipative systems, the frequencies
are expected to depend continuously on the parameters. This means that, for some values of the parameters
the motions can be periodic, and still invariant tori exist. We mention that existence of quasi-periodic solutions
in both conservative and dissipative systems, without assuming a non-resonance condition on the frequencies,
was proved in different contexts~\cite{BC1,BC2,BZ,GV}, but the results do not apply to the class of
systems we are considering.

Of course, without discussing the issue of convergence of the iteration scheme, we cannot conclude that
the solution we look for really exists. Indeed, the solution is defined as the limit --- if it exists ---
of the successive approximations found along the iteration. When computing
the quasi-periodic solution corresponding to the 3:2 resonance, we can claim only
\textit{a posteriori}, after comparing with the results of the numerical investigation, that such a solution really exists. 
The main advantages of the analytical approach, with respect the numerical results, are that:
\begin{enumerate}
\itemsep0em
\item we obtain an explicit formula that approximates reasonably well the quasi-periodic attracting solution;
\item we provide some insight into the nature of such a solution and on the mechanism of its creation.
\end{enumerate}

The rest of the paper is organised as follows. In Section~\ref{sec:2} we define the model, and
give the explicit expressions of the functions $F$ and $G$ appearing in~\eqref{eq:1.1}.
In Section~\ref{sec:3} we study analytically the existence of periodic solutions and find 
explicit expressions for them. These match very well the solutions that are found numerically, represented
in Figure~\ref{figure:1}.
In Section~\ref{sec:5} we study the quasi-periodic attractor close to the resonance 3:2 (the most interesting one,
since it corresponds to the dominant attractor), by applying perturbation theory after
a preliminary step which redefines the unperturbed system: once more we find
good agreement with the numerical results as given in Figure~\ref{figure:2}.
In Section~\ref{sec:6} we discuss the stability of the periodic solutions and the appearance of the quasi-periodic solutions
by Hopf bifurcation, by presenting numerical results which provide further support to the analytical ones. 
In Section~\ref{sec:4} we study analytically the dynamics far from the resonances; we
show that several simplifications can be made, which we justify \textit{a posteriori} by comparison with the
numerical results: the approximate analytical solutions fit the numerical ones closely and reveal a slow, almost linear
decreasing of the velocity.
In Section \ref{sec:8}, inspired by the results of the previous sections, we revisit the problem of determining
the probability of capture of the existing attractors; we find that, assuming an originally faster rotating Mercury,
the probability of capture in the $3$:$2$ resonance is higher than $50\%$.
Finally, Section \ref{sec:7} is devoted to the conclusions and a discussion of open problems.

\section{The spin-orbit model with a realistic tidal torque} \label{sec:2}
\setcounter{equation}{0}

The spin-orbit model describes an asymmetric ellipsoidal celestial body
which moves in a Keplerian elliptic orbit around a central body and rotates around an axis
orthogonal to the orbit plane \cite{D,GP,MD}.

The ordinary differential equation governing the dynamics of the system is
\begin{equation} \label{eq:2.1}
\ddot  \theta = -\zeta \, G( \theta,t) - \eta F(\dot \theta) ,
\end{equation}
where $\theta\in\TTT=\RRR/2\pi\ZZZ$ denotes the angle between the longest axis of the body and the line of apsides,
the parameters $\zeta,\eta\in\RRR_{+}$ are small and the dots denote derivatives with respect to the time $t$.
On the right hand side of \eqref{eq:2.1},
the term with $G(\theta,t)$, that we call the \emph{force} in what follows,
represents the triaxial torque acting on the system,
while the term with $F(\dot \theta)$ models the dissipation due to the tidal torque.
In the literature, $\ddot\theta^{\rm (TRI)}:=-\zeta \, G( \theta,t)$ and
$\ddot\theta^{\rm (TIDE)}:=- \eta F(\dot \theta)$ sometimes are referred to, slightly improperly, as
the triaxial acceleration and the tidal acceleration, respectively.
In this paper we focus on the Mercury-Sun system, even though the ideas could be applied to any system formed
by a satellite orbiting its primary --- of course the tidal model to use strongly depends on the system one is interested in, 
as forcefully pointed out in the recent literature %
\cite{EL,WE,E1,E2,MBE,EM,ME,FM1,FM2,FQEG,MFD}.

The function $G(\theta,t)$ has the form \cite{D,GP,MD}
\begin{equation} \label{eq:2.2}
G(\theta,t) = \sum_{k\in\KK} A_{k} \sin(2\theta- k n t) , \qquad \KK=\{-2,-1,0,1,2,3,4,5,6,7,8\} , 
\end{equation}
where $n$ is the forcing frequency and the coefficients $A_k$, which depend on the eccentricity,
are related to the Hansen coefficients \cite{MD} $G_{20q}$ by $A_{k}=G_{20q}$, with $k=q+2$.
For the Mercury-Sun system, for which $n=26.0879 \, \hbox{ yr}^{-1}$,
the coefficients $A_{k}$ for $k\neq 0$ are given in Table \ref{table:1}, with $A_0=0$.

\begin{table}[ht]
\centering
\begin{center} 
\setlength\tabcolsep{5pt}
\vskip.5truecm
\vrule
\begin{tabular}{|l|l|l|l|l|l} \hline
$k$ & -2  & -1  & 1 & 2 &  3 \\ \hline
$A_{k}$ &  $\st 7.673 \times 10^{-5}$ & $\st 1.865 \times 10^{-4}$ &
$\st -1.023 \times 10^{-1}$ & $\st 8.958 \times 10^{-1}$ & $\st 6.542 \times 10^{-1}$ \\ \hline\hline
$k$ & 4 & 5 & 6 & 7 & 8 \\ \hline
$A_{k}$ & $\st 3.260 \times 10^{-1}$ & $\st 1.380 \times 10^{-1}$ & $\st 5.325 \times 10^{-2}$ &
$\st 1.937 \times 10^{-2}$ & $\st 6.763 \times 10^{-3}$  \\ \hline
\end{tabular}
\hspace{-0.07cm}\vrule
\caption{Values of the coefficients $A_{k}$, with $k\in\KK \setminus\{0\}$, in the case of Mercury-Sun.}
\label{table:1}
\end{center}
\end{table}

The function $F(\dot \theta)$ is given by \cite{EL,E2,M,MBE,NFME}
\begin{equation} \label{eq:2.4}
F(\dot \theta) : = \sum_{k\in\QQ} A_{k}^{2} \, \Xi(\Omega_k) , \qquad \QQ=\{1,2,3,4,5,6,7,8,9\} ,
\end{equation}
where
\begin{equation} \label{eq:2.5}
\Omega_{k}:= n k - 2 \dot \theta , \qquad
\Xi(\omega) = {\rm sgn}(\omega) \, \frac{\mathcal{I}(|\omega|) \, 
|\omega|}{(\mathcal{R}(|\omega|) + \mathcal{A} |\omega|)^2 + \mathcal{I}^2(|\omega|)} ,
\end{equation}
with the functions $\mathcal{I}(\omega)$ and $\mathcal{R}(\omega)$ given by
\begin{subequations} \label{eq:2.6}
\begin{align}
\mathcal{I}(\omega) & = - \tau_M^{-1} - \omega^{1-\alpha} \tau_A^{-\alpha} \sin(\alpha\pi/2) \Gamma(\alpha+1) ,
\label{eq:2.18a} \\
\mathcal{R}(\omega) & = \omega+\omega^{1-\alpha} \tau_A^{-\alpha} \cos(\alpha\pi/2)\,\Gamma(\alpha+1) ,
\label{eq:2.18b}
\end{align}
\end{subequations}
where $\Gamma$ is the Gamma function. The values of the constants appearing in \eqref{eq:2.5} and \eqref{eq:2.6} are
$\alpha=0.2$, $\tau_{M}=\tau_A=500$ and
$\mathcal{A} = 38\pi \mu R^4/3 G M^2 = 15.51726$,
where $\mu=7.967 \times 10^{28}$ km$^{-1}$yr$^{-2}$,
$R=2.44 \times 10^3$ km, $M=3.301 \times 10^{23}$ kg and $G=6.646 \times 10^{-5}$ kg$^{-1}$km$^3$yr$^{-2}$
are the unrelaxed rigidity, the radius of Mercury, the mass of Mercury and the gravitational constant, respectively.
The parameters $\tau_M$ and $\tau_A$ are the Maxwell and Andrade times.
The tidal acceleration $\ddot\theta^{\rm (TIDE)}=-\eta \, F(\dot \theta)$ is plotted, on two different scales,
in \cite[Figures 3 and 4]{BDG2}: there are five pronounced kinks where the function changes sign, the three most significant
being at $\dot\theta/n \approx 1,3/2,2$.

Finally, the parameters $\zeta$ and $\eta$ in \eqref{eq:2.1} are
\begin{subequations} \label{eq:2.3}
\begin{align}
\zeta & := \frac{3}{2} \frac{B - A}{C} n^2 = 0.09545 \hbox{ yr}^{-2}, 
\label{eq:2.3a} \\
\eta & = 0.03096 \hbox{ yr}^{-2}, 
\label{eq:2.3b}
\end{align}
\end{subequations}
where $A<B<C$ are the moments of inertia with respect to the $x$, $y$ and $z$ axes of Mercury.
We refer to \cite{NFME,BDG2} and the references therein for further details of  the spin-orbit model.
All the values of the parameters are taken from~\cite{NFME}.

\section{Method of successive approximations for the periodic attractors} \label{sec:3}
\setcounter{equation}{0}

We set $\zeta=\e$ and $\eta=\gamma\e$ in \eqref{eq:2.1}, with $\gamma = 0.3243$, and write \eqref{eq:2.1} as
\begin{equation} \label{eq:3.0}
\ddot  \theta = -\e \, G( \theta,t) - \e \, \gamma \, F(\dot \theta) ,
\end{equation}
where there is only one parameter. For $\e=0$ the equation is trivially solved.

In this section we study the existence of periodic solutions to \eqref{eq:3.0}, by treating $\e$ as a perturbation
parameter. Since the tidal torque is only $C^{1}$, instead of the usual perturbation theory based on
power series expansions in $\e$, we shall rely on a Picard-like iteration method \cite{CL}, more suited for
differential equations with $C^1$ vector fields, to find successive approximations to the periodic solutions. 

\subsection{Zeroth approximation} \label{subsec:3.1}

We look for a periodic solution which continues the unperturbed one with frequency $\omega_0=p/q \in \QQQ$,
that is, a solution which reduces to $\theta_0(t):=\bar \theta_0 + n \omega_0 t$ as $\e \to 0$. 
Let us write \eqref{eq:3.0} as an integral equation:
\begin{equation} \nonumber 
\theta(t) = \bar  \theta + n \omega_0 t + \bar y t - \e \int_{0}^{t} {\rm d}\tau \int_{0}^{\tau} {\rm d}\tau' 
\left[ G( \theta(\tau'), \tau') + \gamma \, F(\dot \theta(\tau')) \right] ,
\end{equation}
where the constants $\bar \theta$ and $\bar y$ have to be fixed by requiring the solution $\theta(t)$
to be periodic with period $2\pi q/n$.

Let us consider as the zeroth approximation the solution $\theta_0(t)$ of the equation obtained by setting $\e=0$,
$\theta_0(t) = \bar  \theta_0 + n \omega_0 t + \bar y_0 t$,
which yields $\bar y_0=0$, so that
\begin{equation} \label{eq:3.2}
\theta_0(t) = \bar \theta_0 +n \omega_0 t ,
\end{equation}
with $\bar \theta_0$ being arbitrary. This is the unperturbed solution. To fix the value of $\bar \theta_0$ we have
to consider the equation for the first approximation $\theta_1(t)$, which is the solution to the integral equation
\begin{eqnarray} 
\theta_1(t) 
& \!\!\!\! = \!\!\!\! & 
\bar  \theta_1 + n \omega_0 t + \bar y_1 t - \e \int_{0}^{t} {\rm d}\tau \int_{0}^{\tau} {\rm d}\tau' 
\left[ G( \theta_0(\tau'), \tau') + \gamma \, F(\dot \theta_0(\tau')) \right] \nonumber \\
& \!\!\!\! = \!\!\!\! & 
\bar  \theta_1 + n \omega_0 t + \bar y_1 t - \e \int_{0}^{t} {\rm d}\tau \int_{0}^{\tau} {\rm d}\tau' 
\left[ G( \bar \theta_0+ n \omega_0 \tau', \tau') + \gamma \, F(n \omega_0) \right] . \nonumber
\end{eqnarray}
For any periodic function $f$ with period $2\pi q/n$ we denote its average by
\begin{equation} \nonumber
\langle f \rangle := \frac{n}{2\pi q} \int_{0}^{2\pi q/n} {\rm d}t \, f(t) .
\end{equation}
By defining the functions $u_1(t)$, $U_1(t)$ and $\UU_1(t)$ as
\begin{equation} \nonumber
u_1(t) := G( \bar \theta_0 + n \omega_0 t, t) + \gamma \, F(n\omega_0) , \qquad U_1'(t) := u_1(t) , \qquad \UU_1'(t) := U_1(t) ,
\end{equation}
one needs $\langle u_1 \rangle=0$ for the function
\begin{equation} \label{eq:3.4}
\int_{0}^{\tau} {\rm d}\tau' \left[ G( \bar \theta_0 + n \omega_0 \tau', \tau') + \gamma \, F(n \omega_0) \right] = U_1(\tau) - U_1(0)
\end{equation}
to be periodic (and not to grow linearly with $\tau$). This leads to the requirement that
\begin{equation} \label{eq:3.5}
\gamma \, F(n\omega_0) + \frac{n}{2\pi q} \int_{0}^{2\pi q/n} {\rm d}t \, G(\bar \theta_0 + n\omega_0 t, t) = 0 ,
\end{equation}
which, for $G$ as in \eqref{eq:2.2}, yields $2p/q\in\KK$ and
\begin{equation} \label{eq:3.6}
A_{k_0} \sin 2 \bar \theta_0 = - \gamma \, F(n\omega_0) ,\qquad k_0 = \frac{2p}{q} .
\end{equation}
Note that \eqref{eq:3.6} can be satisfied if and only if (i) $q=2$ and $p\in\KK$ and (ii) $\gamma|F(n\omega_0)/A_{k_0}| < 1$.
If these conditions are fulfilled, the constant $\bar \theta_0$ is fixed to one of the $4$ values
\begin{equation} \label{eq:3.7}
\bar \theta_0^{(1)} = \frac{1}{2} {\rm arcsin} \Biggl( - \frac{\gamma \, F(n\omega_0)}{A_{k_0}} \Biggr)
\in \left( -\frac{\pi}{4},\frac{\pi}{4} \right) ,
\quad \bar \theta_0^{(2)} = \frac{\pi}{2} - \bar \theta_0^{(1)} , \quad 
\bar \theta_0^{(3)}=\bar \theta^{(1)} - \pi , \quad \bar \theta_0^{(4)}=\bar \theta_0^{(2)}-\pi ,
\end{equation}
which, at least for small values of $\e$, correspond to two stable and two unstable solutions
(attractors and repellers, respectively); this is a consequence of the Poincar\'e-Birkhoff theorem \cite{AKN,C}.
We note at this point that bifurcation phenomena may occur when increasing the value of $\e$
(see Section~\ref{sec:6} below).

Therefore the zeroth approximation is given by~\eqref{eq:3.2},  with $\bar \theta_0$ given by one of
the four values in~\eqref{eq:3.7}.

\subsection{First approximation} \label{subsec:3.2}

Now we want to compute the first approximation $\theta_1(t)$. Once $\bar \theta_0$ has been fixed in such
a way that \eqref{eq:3.5} is satisfied, we can compute $U_1$ in \eqref{eq:3.4}.
If we require $\langle U_1 \rangle$ to vanish (for the function $\UU_1(t)$ to be bounded and hence periodic),
we obtain, for $G$ as in \eqref{eq:2.2},
\begin{equation} \nonumber 
U_1(t) = - \sum_{k\in\KK_0} \frac{A_k}{(2\omega_0-k)n} \cos(2 \bar \theta_0+(2\omega_0-k)nt) ,
\end{equation}
where $\KK_0:=\KK\setminus\{k_0\}$. Then, by fixing $\bar y_1 = - U_1(0)$,
and choosing $\UU_1(t)$ so that its average $\langle \UU_1 \rangle$ vanishes, one finds
$\theta_1(t) = \Th_1 + n \omega_0 t - \e \UU_1(t)$, where $\Th_1 := \bar \theta_1 - \e \UU_1(0)$
and, for $G$ as in \eqref{eq:2.2},
\begin{equation} \nonumber
\UU_1(t) = - \sum_{k\in\KK_0} \frac{A_k}{(2\omega_0-k)^2n^2} \sin(2 \bar \theta_0+(2\omega_0-k)nt) .
\end{equation}
Thus, we obtain
\begin{equation} \label{eq:3.8}
\theta_1(t) = \Th_1 + n\omega_0 t + \e \sum_{k\in\KK_0} \frac{A_k}{(2\omega_0-k)^2n^2} \sin(2 \bar \theta_0+(2\omega_0-k)nt) ,
\end{equation}
with $\Th_1$ to be determined by requiring that the second order $\theta_2(t)$ has period $2\pi q/n$ as well.

The second approximation $\theta_2(t)$ is the solution to the equation
\begin{equation} \label{eq:3.9}
\theta_2(t) = \bar  \theta_2 + n \omega_0 t + \bar y_2 t - \e \int_{0}^{t} {\rm d}\tau \int_{0}^{\tau} {\rm d}\tau' 
\left[ G( \theta_1(\tau'), \tau') + \gamma \, F(\dot \theta_1(\tau')) \right] ,
\end{equation}
where $\theta_1(t)$ is the solution \eqref{eq:3.8} found at the first iterative step.
Therefore, for the solution to be periodic, we need
\begin{equation} \label{eq:3.10}
\frac{n}{2\pi q} \int_{0}^{2\pi q/n} {\rm d}\tau'
\left[ G( \Th_{1}+n\omega_0 \tau' + \xi_1(\tau'),\tau') + \gamma \, F(n\omega_0 + \dot \xi_1(\tau')) \right] =  0 ,
\end{equation}
where $\xi_1(t) := - \e  \UU_1(t)$. If we rewrite \eqref{eq:3.10} as
\begin{equation} \label{eq:3.11}
\mathcal{G}(\Th_1) := \frac{n}{2\pi q} \int_{0}^{2\pi q/n} {\rm d}\tau'
\, G( \Th_{1}+n\omega_0 \tau' + \xi_1(\tau'),\tau') =
- \frac{n}{2\pi q} \int_{0}^{2\pi q/n} {\rm d}\tau' \, \gamma \, F(n\omega_0 + \dot \xi_1(\tau')) ,
\end{equation}
we see we have to invert the function $\mathcal{G}$ to find $\Th_1$.
By expanding and using the fact that $\xi_1=O(\e)$, we find
\begin{equation} \nonumber
G( \Th_{1}+n\omega_0 \tau' + \xi_1(\tau'),\tau') =
G( \Th_{1}+n\omega_0 \tau' ,\tau') + \Gamma(\tau',\e), \quad
\Gamma(\tau',\e) = \partial_\theta G( \Th_{1}+n\omega_0 \tau' ,\tau') \, \xi_1(\tau')  + O(\e^2) ,
\end{equation}
so that we can write
\begin{equation} \label{eq:3.12}
\mathcal{G}(\Th_1) =
\frac{n}{2\pi q} \int_{0}^{2\pi q/n} {\rm d}\tau' G( \Th_{1}+n\omega_0 \tau',\tau') + O(\e) =
A_{k_0} \sin 2 \Th_1 + O(\e) .
\end{equation}
By neglecting the corrections of order $\e$ in \eqref{eq:3.12} and defining
\begin{equation} \nonumber
J(n\omega_0) := \frac{n}{2\pi q} \int_{0}^{2\pi q/n} {\rm d}\tau' \, \gamma \, F(n\omega_0 + \dot \xi_1(\tau')) ,
\end{equation}
we find
\begin{equation} \label{eq:3.13}
A_{k_0} \sin 2 \Th_1 = -J(n\omega_0) ,\qquad k_0 = \frac{2p}{q} ,
\end{equation}
where again one must have $q=2$ and $p\in\KK$, so that 
\begin{equation} \label{eq:3.14}
\Th_1^{(1)} = \frac{1}{2} {\rm arcsin} \Biggl( - \frac{J(n\omega_0)}{A_{k_0}} \Biggr) \in \left( -\frac{\pi}{4},\frac{\pi}{4} \right) ,
\quad \Th_1^{(2)} = \frac{\pi}{2} - \Th_1^{(1)} , \quad 
\Th_1^{(3)}=\tilde \theta^{(1)} - \pi , \quad \Th_1^{(4)}=\Th_1^{(2)}-\pi .
\end{equation}
The values of the integrals $J(n\omega_0)$ in \eqref{eq:3.13} are given in Table \ref{table:2}, where they are also
compared to the values $\gamma F(n \omega_0)$. It is evident that there are appreciable discrepancies between the two
values for $\omega=1$, 3/2 and 2, which correspond to the major kinks of the function $F$.

\begin{table}[H]
\centering
\begin{center}
\setlength\tabcolsep{5pt}
\vskip.4truecm
\vrule
\begin{tabular}{|l|l|l|l|l|l} \hline
$\omega_0$ & -1  & -1/2  & 1/2 & 1 & 3/2  \\ \hline
$\gamma \, F(n\omega_0)$ & $\st -5.09302 \times 10^{-5}$ & $\st -5.36152 \times 10^{-5}$ &
$\st -6.387 \times 10^{-5}$ & $\st -2.639 \times 10^{-5}$ & $\st 3.429\times 10^{-5}$ \\ \hline
$J(n\omega_0)$ & $\st -5.09277 \times 10^{-5}$ & $\st -5.36126 \times 10^{-5}$ &
$\st -6.073 \times 10^{-5}$ & $\st 1.390 \times 10^{-4}$ & $\st 1.163 \times 10^{-4}$ \\ \hline\hline
$\omega_0$ &  2 & 5/2 & 3 & 7/2 & 4 \\ \hline
$\gamma \, F(n\omega_0)$ & $\st 5.557 \times 10^{-5}$ & $\st 5.646 \times 10^{-5}$ &
$\st 5.363 \times 10^{-5}$  & $\st 5.102 \times 10^{-5}$  & $\st 4.898\times 10^{-5}$  \\ \hline
$J(n\omega_0)$ & $\st -2.143 \times 10^{-5}$ & $\st 5.172 \times 10^{-5}$ &
$\st 5.332 \times 10^{-5}$  & $\st 5.100 \times 10^{-5}$  & $\st 4.898 \times 10^{-5}$  \\ \hline
\end{tabular}
\hspace{-0.07cm}\vrule
\caption{Values of $\gamma \, F(n\omega_0)$ and $J(n\omega_0)$.
More than four significant figures are needed only for the retrograde resonances, 
in order to observe in practice the difference
between the corresponding curves.}
\label{table:2}
\end{center}
\end{table}

The values of the constants $\Th_1$ in the interval $(-\pi/4,\pi/4)$ are given in  Table \ref{table:3}, where
the corresponding values $\bar \theta_0$ are also given. Once more the difference between the
two values $\bar \theta_0$ and $\Th_1$ is larger for resonances such as $\omega_0=1$, $\omega_0=3/2$
and $\omega_0=2$, where the kinks of the function $F$ are more pronounced.

\begin{table}[ht]
\centering
\vskip.4truecm
\begin{center}
\setlength\tabcolsep{5pt}
\vrule
\begin{tabular}{|l|l|l|l|l|l} \hline
$\omega_0$ & -1  & -1/2  & 1/2 & 1 & 3/2  \\ \hline
\raisebox{2.2ex}{}
$\bar \theta_0$ & $\st 3.62911 \times 10^{-1}$ & $\st 1.45808 \times 10^{-1}$ &
$\st -3.123 \times 10^{-4}$ & $\st 1.473 \times 10^{-5}$ & $\st -2.621 \times 10^{-5}$ \\ \hline
\raisebox{2.2ex}{}
$\Th_1$ & $\st 3.62889 \times 10^{-1}$ & $\st 1.45801 \times 10^{-1}$ &
$\st  -2.969 \times 10^{-4}$ & $\st -7.758 \times 10^{-5}$ & $\st -8.888 \times 10^{-5}$ \\ \hline\hline
$\omega_0$ &  2 & 5/2 & 3 & 7/2 & 4 \\ \hline
\raisebox{2.2ex}{}
${\bar \theta_0}$ & $\st - 8.541 \times 10^{-5}$ & $\st - 2.046 \times 10^{-4}$ &
$\st - 5.035\times 10^{-4}$ & $\st - 1.317 \times 10^{-3}$ & $\st -3.621 \times 10^{-3}$ \\ \hline
\raisebox{2.2ex}{}
$\Th_1$ & $\st 3.259 \times 10^{-5}$ & $\st -1.874 \times 10^{-4}$ &
$\st -5.006 \times 10^{-4}$  & $\st -1.316 \times 10^{-3}$  & $\st -3.621 \times 10^{-3}$  \\ \hline
\end{tabular}
\hspace{-0.07cm}\vrule
\caption{Values of the constants $\bar \theta_0$ and $\Th_1$.
More than four significant figures are needed only for the retrograde resonances, in order to observe in practice the difference
between the corresponding curves.}
\label{table:3}
\end{center}
\end{table}

\begin{rmk} \label{rmk:1}
\emph{
One may wonder why in \eqref{eq:3.11} we Taylor-expand the function $G$ but not the function $F$. This is due to the
fact that $G(\theta,t)$ is a smooth function (in fact it is analytic), while $F$ is only $C^1$. Moreover,
the first derivative of $F$ is very large at some resonances such as $\om_0=1$ and $\om_0=3/2$;
an explicit computation gives $\eta\partial_{\dot\theta}F(n) = 24.8421$ and
$\eta\partial_{\dot\theta}F(3n/2) = 13.2493$. 
Thus, since the size of $\xi_1(t)$ is of order $\e$, for $\e=\zeta$ fixed as in \eqref{eq:2.3},
the function $G(\Th_{1}+n\omega_0 t + \xi_1(t),t)$ is well approximated by $G( \Th_{1}+n\omega_0 t,t)$,
while $F(n\omega_0 + \dot \xi_1(\tau'))$ can be appreciably different from $F(n\omega_0)$.
}
\end{rmk}

\begin{rmk} \label{rmk:2}
\emph{
The reason why we neglect the corrections of order $\e$ in \eqref{eq:3.12} is that, then,
the implicit function equation \eqref{eq:3.11}
is easily solved. Of course, in so doing, an error is introduced. However, we argue that such an error is of the same
order of magnitude of the terms we are disregarding by stopping Picard's iteration at the first step.
Indeed, as the difference between $\theta_1(t)$ and $\theta_0(t)$ is of size $\e$,
so the difference between $\theta_2(t)$ and $\theta_1(t)$ is expected to be of size $\e^2$.
Taking into account the terms $O(\e)$ in \eqref{eq:3.12} would lead to values
of $\Th_1$ which differ from those given by \eqref{eq:3.14} by terms of size $\e^2$ --- the same size
as the corrections to be found in the next iterative step.
}
\end{rmk}

\begin{figure}[!ht]
\centering 
\vspace{-0.2cm}
\includegraphics[width=6.2in]{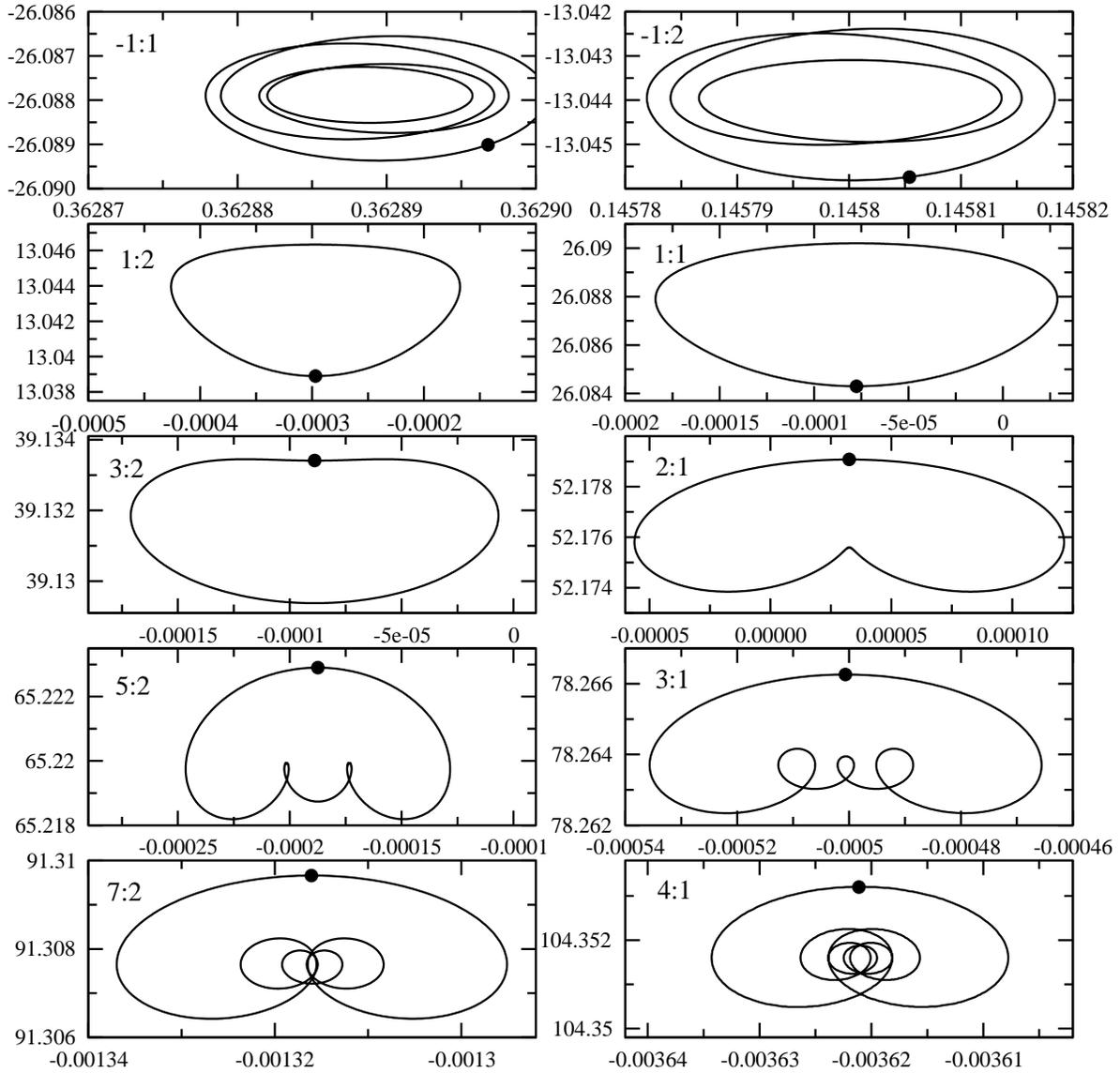}
\caption{Periodic solutions according to the first approximations \eqref{eq:3.8} with the constants
$\bar\theta_0$ and $\Th_1$ taken from Table \ref{table:3}; in each figure the dot represents the initial condition.
The arrangement of the resonances and the range of the variables are both as in Figure \ref{figure:1}.}
\label{figure:4}
\end{figure} 


In conclusion, the first approximation $\theta_1(t)$ is given by \eqref{eq:3.8}, with the constants $\bar\theta_0$ and
$\Th_1$ taken from Table \ref{table:3}. By defining the libration  $z(t):=\theta_1(t)-n\omega_0 t$, one finds that,
in the $(z,\dot \theta)$-plane, the periodic attractors are described by the curves in Figure \ref{figure:4}.
There is an attractor for any value of $\om_0=p/q$, with $q=2$ and $p\in\KK\setminus\{0\}$.
Each attractor is characterised by the property that it describes an oscillation $\theta(t)$ with $\dot\theta(t)$ close
to $n\om_0$, while $\theta(t) - n \om_0 t$ moves around $\Th_{1}$. Note that the amplitude of the libration
is essentially determined by the triaxial torque, as expected for relatively cold celestial bodies with high viscosity and
large Maxwell time $\tau_M$ \cite{MFD}; however its centre is fixed by the tidal torque, through \eqref{eq:3.13}.

\begin{rmk} \label{rmk:3}
\emph{
The first approximation is not a first order perturbation theory solution.
Indeed, if one considered the Taylor expansion in $\e$
of the solution and kept only the terms up to the first order, then one should replace $J(n\omega)$ with $\gamma \, F(n\omega)$
in \eqref{eq:3.11}, so that $\Th_1$ would become a correction of order $\e$ of $\bar \theta_0$.
Thus, the solution would be still of the form \eqref{eq:3.8},
but with a different value for $\Th_1$. An explicit computation shows that the value inside the interval
$(-\pi/4,\pi/4)$ would be appreciably smaller than $\Th_1$, as given in Table \ref{table:3},
so that the curves would appear slightly shifted in the $z$-direction with respect to those represented in Figure \ref{figure:4}.
This effect would be noticeable in all periodic attractors.
}
\end{rmk}

\subsection{Comparison with the numerical results} \label{sec:3.3}
\setcounter{equation}{0}

For all the numerical computations --- both in this section and in the following ones --- we rely on the
fast numerical integrator introduced in \cite{BDG1,BDG2}, to which we refer for details. In order to compute the attractors,
we follow the time evolution of a large number of initial conditions in phase space (more than 50 000)
so as to determine the asymptotic behaviour of the corresponding trajectories --- see also \cite{BBDGG}.
By using a large number of initial conditions we argue that even if, besides those that have been detected, 
other attractors existed, they would be irrelevant, since they would attract only a negligible fraction of trajectories.
Of course, if other periodic attractors exist, they are not obtained by continuation from the unperturbed
solutions of the form \eqref{eq:3.2}. The periodic solutions which are found numerically are represented
in Figure \ref{figure:1}. Not all of them are attractors (see Section \ref{sec:6} below):
the periodic solutions which are unstable may be detected by a fixed point method by considering the
corresponding Poincar\'e section.

A comparison with the analytical results shows that there is a very good agreement
so that we conclude that the first order approximation $\theta_1(t)$ provides an accurate description of the solution.

Actually, a shift is observed in the case of the resonances $\omega_0=-1/2$ and, particularly, $\omega_0=-1$.
The reason behind that is very likely due to
the smallness of the coefficients $A_k$ corresponding to the retrograde resonances ($\omega_0=-1/2$ and $\omega_0=-1$),
which makes the dependence of $\Th_1$ very sensitive to the exact value of the integral $J(n\omega_0)$.
For instance an error in the fourth decimal digit, which has no effect for the prograde resonances,
is able to produce an appreciable shift for the retrograde resonances.

\section{The quasi-periodic attractor close to the 3:2 resonance} \label{sec:5}
\setcounter{equation}{0}

We study now the existence of quasi-periodic solutions to \eqref{eq:2.1} describing invariant tori.
Consider the ordinary differential equation \eqref{eq:2.1},
with the functions  $G(\theta,t)$ and $F(\dot \theta)$ as in \eqref{eq:2.2} and \eqref{eq:2.4},
respectively. We look for a solution $\theta(t)$ with $\dot\theta(t) \approx 3n/2$, which suggests setting
\begin{equation} \label{eq:5.1}
\theta(t) := \frac{3}{2} n t + \frac{1}{2} \xi(t) .
\end{equation}
In terms of $\xi$ the equation of motion becomes
\begin{equation} \nonumber
\ddot  \xi = - 2 \zeta \sum_{k\in\KK} A_{k} \sin( \xi - (k-3) n t) - 2 \eta \, \Phi(\dot \xi) , \qquad
\Phi(\dot \xi) := F \Bigl( \frac{3n + \dot \xi}{2} \Bigr) ,
\end{equation}
which we rewrite as
\begin{equation} \label{eq:5.2}
\ddot  \xi = - \omega^2 \sin \xi - 2 \zeta \sum_{k\in\KK_{0}} A_{k} \sin( \xi - (k-3) n t) - 2 \eta \, \Phi(\dot \xi) ,
\end{equation}
where $\omega := \sqrt{2 \zeta A_{3}} =  0.3534$
and $\KK_{0}=\KK \setminus \{ k_0 \}$, with $k_{0}=3$.

If in \eqref{eq:5.2} we set $\zeta=\e$ and $\eta=\gamma\e$, without affecting $\om$,
\eqref{eq:5.2} becomes
\begin{equation} \label{eq:5.3}
\ddot  \xi = - \omega^2 \sin \xi - 2 \e \sum_{i\in\II_{0}} B_{i} \sin( \xi - in t) - 2 \gamma \, \e \Phi(\dot \xi) ,
\end{equation}
where we have defined $B_{i}:=A_{i+3}=G_{20i+1}$ and $\II_{0}=\{\pm1,\pm2,\pm3,\pm4,\pm5\}$.
We would like to study \eqref{eq:5.3} by considering $\omega$ as a parameter independent of $\e$.
Thus, \eqref{eq:5.3} could be seen as a perturbation of the pendulum equation, $\ddot  \xi = - \omega^2 \sin \xi$,
to which it reduces as $\e=0$. However, since we look for a solution to \eqref{eq:5.3} with $\xi$ close to $0$,
to simplify the analysis we proceed in a slightly different way. Starting from \eqref{eq:5.2}, 
we make a linear approximation of the unperturbed system by splitting $\sin\xi$ into two terms
and re-writing \eqref{eq:5.2} as
\begin{equation} \label{eq:5.4}
\ddot  \xi = - \omega^2 \xi - 2 \e \sum_{i\in\II_{0}} B_{i} \sin( \xi - in t) - 2 \gamma \, \e \Phi(\dot \xi) -
2 A_3 \e \left( \sin \xi - \xi \right) ,
\end{equation}
where again $\om$ is to be considered fixed, while $\e$ can be varied. Of course, eventually, we have to
fix $\e$ to the value such that $2\e A_3=\om^2$. Note that  a similar approach has been proposed in~\cite{W,GCEP}
to study the secondary resonances in non-dissipative spin-orbit models with $\e$ close to a rational number
and small values of eccentricity --- a situation very far from that of Mercury.

In general, quasi-periodic solutions to \eqref{eq:5.4} will have
two frequencies $\omega_L$ and $n$, with $n$ fixed and $\omega_L$ close to $\omega$ --- see \cite{BBDGG,GBD,WBG}.
We call $\om_{L}$ the low frequency and $n$ the high frequency, because $\om_{L}\approx \om$ and $n=73.82\om$.
Hence we write the solution as $\xi(t)=X(\oo_L t)$, where $\oo_L:=(\om_{L},n)$.
To take into account the possible dependence of frequency on the perturbation we write \eqref{eq:5.4} as
\begin{equation} \label{eq:5.5}
\ddot  \xi = - \omega_{L}^2 \xi - \mu \e\,  \xi - 2 \e \sum_{i\in\II_{0}} B_{i} \sin( \xi - in t) - 2 \gamma \, \e \Phi(\dot \xi) 
- 2 A_3 \e \left( \sin \xi - \xi \right) ,
\end{equation}
where $\om^2=\om_{L}^2+\mu\,\e$, with $\om_L$ and $\mu$ constants to be determined. As a further simplification,
we approximate $\sin\xi \approx \xi -\xi^3/6$, so that \eqref{eq:5.5} becomes
\begin{equation} \label{eq:5.6}
\ddot  \xi = - \omega_{L}^2 \xi - \mu \e\,  \xi - 2 \e \sum_{i\in\II_{0}} B_{i} \sin( \xi - in t) - 2 \gamma \, \e \Phi(\dot \xi) + \frac{A_3 \e}{3} \, \xi^3  ,
\end{equation}

We look for approximate solutions to \eqref{eq:5.6} using a Picard-like iteration scheme as in Section \ref{sec:3}.
The strategy is the following. The zeroth approximation is the solution to \eqref{eq:5.6} with $\e=0$.
Suppose that at step $k\ge 0$ we have found an approximation $\xi_{k}(t)=X_{k}(\oo_L t)$, with
frequency vector $\oo_{L}=(\om_L,n)$, for some $\om_L$ close to $\om$. As we shall see such a solution
depends on an arbitrary parameter $C_k$. Then we compute the $(k+1)-$th approximation as the solution of
\begin{equation} \nonumber
\ddot  \xi = - \omega_{L}^2 \xi - \mu \e\,  \xi_{k}(t) - 2 \e \sum_{i\in\II_{0}} B_{i} \sin( \xi_{k}(t) - in t) - 
2 \gamma \, \e \Phi(\dot \xi_{k}(t)) + \frac{A_3 \e}{3} \xi^3_{k}(t) ,
\end{equation}
where we consider $\om_{L}$ as a parameter related to $\mu$ through the relation $\om^2=\om_{L}^2+\mu\e$, with $\om$ given.
For such a solution to be bounded and hence quasi-periodic with frequency vector $\oo_L$,
we have to fix the parameters $C_k$, $\mu$ and $\om_L$, the latter two so as to also satisfy the constraint that
$\om^2=\om_{L}^2+\mu\e$; the value of the low frequency $\om_{L}$
will be determined by an implicit function problem. The corresponding solution will be a quasi-periodic
function $\xi_{k+1}(t)=X_{k+1}(\oo_L t)$, with $\oo_L$ slightly different from the value found at the previous step,
depending on a new arbitrary constant $C_{k+1}$ to be fixed at the next iteration step
to a value close to $C_{k}$.

\subsection{Zeroth approximation} \label{sec:5.1}

The zeroth approximation $\xi_0(t)$ is the solution to the equation
\begin{equation} \label{eq:5.7}
\ddot  \xi = - \omega_L^2 \xi ,
\end{equation}
obtained from \eqref{eq:5.6} by setting $\e=0$. The solution to \eqref{eq:5.7} is
%
$\xi_{0}(t)=\bar\xi_0 \cos \om_{L} t + \bar y_0 \sin \om_{L} t $,
%
where $\bar\xi_0$ and $\bar y_{0}$ are related to the initial data through the relation $\bar\xi_0=\xi(0)$
and $\bar y_{0}=\dot \xi(0)/\om$. 

For any initial datum, we can write
\begin{equation} \label{eq:5.9}
\xi_{0}(t)= C_0 \, \sin(\om_{L} t + \varphi_0) ,
\end{equation}
where $C_0>0$ and $\varphi_0$ is the initial phase.
By varying the initial phase $\varphi_0$ the corresponding trajectories in the space $(\xi,y,t)$
describe an invariant torus which appears as a right circular cylinder with axis along the $t$-axis and radius $C_0$.

For $\e\neq0$ we expect the cylinder to persist, albeit deformed with respect to the unperturbed case.
This allows us to fix arbitrarily the phase $\varphi_0$, which plays no role: we can put it equal to $0$ for convenience.
Then the torus is parameterised in terms of the value $C_0\in\RRR_{+}$ at which it crosses the
positive $\xi$-axis (corresponding to the phase $\varphi_0=0$).
We expect that the values of $C_0$ and $\om_{L}$, which are arbitrary for $\e=0$, will be fixed when we set $\e\neq0$
and take into account the dissipation by the requirement that the perturbed solution is still bounded.

\begin{rmk} \label{rmk:4}
\emph{
Since we are making a linear approximation, all unperturbed solutions have the same frequency $\om_{L}$.
Therefore all the approximate quasi-periodic solutions \eqref{eq:5.9}, by construction, will have 
the same frequency vector $(\om_{L},n)$. As we shall see, if we wish to compute the higher order approximations,
we need to take into account the change of frequency $\om_{L}$ with respect
to the linearised unperturbed system --- see \cite{BG,BDG1}. In this regard, we note that,
even though usually perturbation theory computations are easier in terms of action-angle variables,
in our case it is more convenient to work with Cartesian coordinates, because we are looking
for a solution around the origin, where the action-angle variables are singular --- see \cite{CGP}
and references therein for similar comments.
}
\end{rmk}

\subsection{First approximation} \label{sec:5.2}

The first order approximation is obtained as the bounded solution to the equation
\begin{equation} \label{eq:5.10}
\ddot  \xi = - \omega_{L}^2 \xi - \mu \e\,  \xi_{0}(t) - 2 \e \sum_{i\in\II_{0}} B_{i} \sin( \xi_{0}(t) - in t) - 
2 \gamma \, \e \Phi(\dot \xi_{0}(t)) + \frac{A_3 \e}{3} \, \xi^3_{0}(t) ,
\end{equation}
with $\xi_0(t)$ given by \eqref{eq:5.9}.
The nonhomogeneous linear equation \eqref{eq:5.10} can be written as a
first order differential equation in $\RRR^2$,
\begin{equation} \label{eq:5.11}
\begin{cases}
\dot \xi = \omega_L y , & \\
\displaystyle{\dot y = - \om_L \xi - \frac{\mu \e}{\om_{L}} \, \xi_{0}(t) + \frac{2\e}{\om_{L}}
\sum_{i\in\II_{0}} B_{i} \sin (in t - \xi_{0}(t)) - \frac{2 \gamma \, \e}{\om_{L}} \Phi(\dot\xi_{0}(t)) } 
+ \frac{A_3 \e}{3\om_L} \, \xi^3_{0}(t) . & \end{cases}
\end{equation}
where $\xi_0(t)=C_0\sin\om_{L} t$ and $\dot\xi_0(t)=C_0\om_{L} \cos \om_{L} t$.
More generally, the approximation at step $k$ is defined as the solution to
\begin{equation} \label{eq:5.10bis}
\begin{cases}
\dot  \xi =\omega_L y , & \\
\displaystyle{\dot y = - \omega_{L}^2 \xi - \mu \e\,  \xi_{k-1}(t) - 2 \e \sum_{i\in\II_{0}} B_{i} \sin( \xi_{k-1}(t) - in t) - 
2 \gamma \, \e \Phi(\dot \xi_{k-1}(t)) + \frac{A_3 \e}{3} \, \xi^3_{k-1}(t) } ,
\end{cases}
\end{equation}
where $\xi_{k-1}(t)$ is the approximation found at step $k-1$. At any step $\om_L$ is to be considered
a free parameter, close to $\om$, to be fixed by requiring the solution to be bounded.

The general solution to \eqref{eq:5.10bis} is
\begin{eqnarray}
& & \left( \begin{matrix} \xi_{k}(t) \\ y_{k}(t) \end{matrix} \right) =
\left( \begin{matrix} \cos\om_{L} t & \sin \om_{L} t \\ -\sin \om t & \cos \om_{L} t \end{matrix} \right)
\Biggl[ \left( \begin{matrix} \bar \xi_{k} \\ \bar y_{k} \end{matrix} \right)  -
\frac{2\gamma\,\e}{\om_{L}} \int_{0}^{t} {\rm d}\tau 
\left( \begin{matrix} \cos\om_{L} \tau & - \sin \om_{L} \tau \\ \sin \om_{L} \tau & \cos \om_{L} \tau \end{matrix} \right)
\left( \begin{matrix} 0 \\  \Phi(\dot\xi_{k-1}(\tau)) \end{matrix} \right) 
\nonumber \\
& & \hspace{1.4cm} 
- \frac{\mu\e}{\om_{L}} \int_{0}^{t} {\rm d}\tau 
\left( \begin{matrix} \cos\om_{L} \tau & - \sin \om_{L} \tau \\ \sin \om_{L} \tau & \cos \om_{L} \tau \end{matrix} \right)
\left( \begin{matrix} 0 \\  \xi_{k-1}(\tau) \end{matrix} \right) 
+ \frac{A_3 \e}{3\om_{L}} \int_{0}^{t} {\rm d}\tau 
\left( \begin{matrix} \cos\om_{L} \tau & - \sin \om_{L} \tau \\ \sin \om_{L} \tau & \cos \om_{L} \tau \end{matrix} \right)
\left( \begin{matrix} 0 \\  \xi^3_{k-1}(\tau) \end{matrix} \right) \nonumber \\
& & \hspace{1.4cm} 
+ \frac{2\e}{\om_{L}} \sum_{i\in\II_{0}} B_{i} \int_{0}^{t} {\rm d}\tau 
\left( \begin{matrix} \cos\om_{L} \tau & - \sin \om_{L} \tau \\ \sin \om_{L} \tau & \cos \om_{L} \tau \end{matrix} \right)
\left( \begin{matrix} 0 \\  \sin (i n \tau - \xi_{k-1}(\tau))\end{matrix} \right) \Biggr] , \nonumber
\end{eqnarray}
where $(\bar\xi_k,\bar y_k)$ is the initial condition. This leads to 
\begin{eqnarray} \label{eq:5.11a}
\xi_{k}(t) & \!\!\! = \!\!\! &
\bar\xi_{k} \cos \om_{L} t + \bar y_{k} \sin \om_{L} t \nonumber \\
& \!\!\! - \!\!\! & \frac{2\gamma\e}{\om_{L}} 
\left( \sin \om_{L} t  \!\! \int_{0}^{t} \!\! {\rm d}  \tau \, \cos \om_{L} \tau \; \Phi(\dot\xi_{k-1}(\tau)) -
\cos \om_{L} t  \!\! \int_{0}^{t} \!\! {\rm d} \tau \, \sin \om_{L} \tau \; \Phi(\dot\xi_{k-1}(\tau)) \right) \nonumber \\
& \!\!\! - \!\!\! & \frac{\mu\e}{\om_{L}}  
\left( \sin \om_{L} t  \!\! \int_{0}^{t} \!\! {\rm d}  \tau \, \cos \om_{L} \tau \; \xi_{k-1}(\tau) -
\cos \om_{L} t  \!\! \int_{0}^{t} \!\! {\rm d} \tau \, \sin \om_{L} \tau \; \xi_{k-1}(\tau) \right) \\
& \!\!\! + \!\!\! & \frac{A_3 \e}{3 \om_{L}}  
\left( \sin \om_{L} t  \!\! \int_{0}^{t} \!\! {\rm d}  \tau \, \cos \om_{L} \tau \; \xi^3_{k-1}(\tau) -
\cos \om_{L} t  \!\! \int_{0}^{t} \!\! {\rm d} \tau \, \sin \om_{L} \tau \; \sin \xi^3_{k-1}(\tau) \right) \nonumber \\
& \!\!\! + \!\!\! & \frac{2\e}{\om_{L}}  \sum_{i\in\II_{0}} \! B_{i} \!
\left( \sin \om_{L} t  \!\! \int_{0}^{t} \!\! {\rm d}  \tau \, \cos \om_{L} \tau \, \sin (in\tau-\xi_{k-1}(\tau)) -
\cos \om_{L} t  \!\! \int_{0}^{t} \!\! {\rm d} \tau \, \sin \om_{L} \tau \, \sin (in\tau-\xi_{k-1}(\tau)) \right) . \nonumber
\end{eqnarray}
For $k=1$ we write $\bar\xi_1 \cos \om_{L} t + \bar y_1 \sin \om_{L} t = C_{1} \, \sin ( \om_{L} t + \varphi_1 )$
and set $\varphi_1=0$ as done for $k=0$. If we define
\begin{subequations} \label{eq:5.11b}
\begin{align} 
M_{k,c}(t) & := 
\int_{0}^{t} \!\! {\rm d}  \tau \, \cos \om_{L} \tau \left[
2 \sum_{i\in\II_{0}} \! B_{i} \sin (in\tau-\xi_{k-1}(\tau)) - 2 \gamma \; \Phi(\dot\xi_{k-1}(\tau)) - 
\mu \xi_{k-1}(\tau) + \frac{A_3}{3} \, \xi^3_{k-1}(\tau) \right] ,
\label{eq:5.11ba} \\
M_{k,s}(t) & :=
\int_{0}^{t} \!\! {\rm d}  \tau \, \sin \om_{L} \tau \left[
2 \sum_{i\in\II_{0}} \! B_{i} \sin (in\tau-\xi_{k-1}(\tau)) - 2 \gamma \; \Phi(\dot\xi_{k-1}(\tau)) - 
\mu \xi_{k-1}(\tau) + \frac{A_3}{3} \, \xi^3_{k-1}(\tau) \right] ,
\label{eq:5.11bb}
\end{align}
\end{subequations}
we can write the solution $\xi_1(t)$ as
\begin{equation} \label{eq:5.12}
\xi_{1}(t) = C_1 \sin \om_{L} t + \Xi_{1}(t) , \qquad
\Xi_{1}(t) := \frac{\e}{\om_{L}} \Bigl( \sin \om_{L} t \,  M_{1,c}(\tau) - \cos \om_{L} t \,  M_{1,s}(\tau) \Bigr) .
\end{equation}

For the functions in~\eqref{eq:5.12} to be bounded one needs
\begin{equation} \label{eq:5.13}
I_{1,c} := \lim_{T\to+\io} \frac{M_{1,c} (T)}{T} = 0 , \qquad I_{1,s} := \lim_{T\to+\io} \frac{M_{1,s} (T)}{T} = 0 .
\end{equation}
The two limits can be computed as
\begin{eqnarray}
I_{1,c} & \!\!\!\! = \!\!\!\! & 
\int_{0}^{2\pi} \frac{{\rm d}\psi_1}{2\pi} \int_{0}^{2\pi} \frac{{\rm d}\psi_2}{2\pi} \; \cos \psi_1 \times \nonumber \\
& & \qquad \times \left[
2 \sum_{i\in\II_{0}} \! B_{i} \sin (i\psi_2-C_0 \sin \psi_1 ) - 2 \gamma \; \Phi(C_0\om_{L} \cos \psi_1 ) - 
\mu \, C_0 \sin \psi_1 + \frac{A_3}{3} C_0^3 \sin^3 \psi_1 \right] ,
\nonumber \\
I_{1,s} & \!\!\!\! = \!\!\!\! & 
\int_{0}^{2\pi} \frac{{\rm d}\psi_1}{2\pi} \int_{0}^{2\pi} \frac{{\rm d}\psi_2}{2\pi} \; \sin\psi_1 \times \nonumber \\
& & \qquad \times \left[
2 \sum_{i\in\II_{0}} \! B_{i} \sin (i\psi_2-C_0 \sin \psi_1 ) - 2 \gamma \; \Phi(C_0\om_{L} \cos \psi_1 ) - 
\mu \, C_0 \sin \psi_1 + \frac{A_3}{3} C_0^3 \sin^3 \psi_1 \right] ,
\nonumber
\end{eqnarray}
since the time average of a quasi-periodic function equals the average over the torus. After expanding
\begin{equation} \nonumber
\sin (i\psi_2-C_0 \sin \psi_1 ) =
\sin i\psi_2 \cos (C_0 \sin \psi_1 ) - \cos i\psi_2 \sin (C_0 \sin \psi_1 )
\end{equation}
and using
\begin{equation} \nonumber
\int_{0}^{2\pi} \frac{{\rm d}\psi_1}{2\pi} \int_{0}^{2\pi} \frac{{\rm d}\psi_2}{2\pi} \; 
\sin \psi_1 \Phi(C_0\om_{L} \cos \psi_1 ) = 0
\end{equation}
by parity, and
\begin{equation} \nonumber
\int_{0}^{2\pi} \frac{{\rm d}\psi_1}{2\pi} \int_{0}^{2\pi} \frac{{\rm d}\psi_2}{2\pi} \; \cos \psi_1 \sin \psi_1 = 
\int_{0}^{2\pi} \frac{{\rm d}\psi_1}{2\pi} \int_{0}^{2\pi} \frac{{\rm d}\psi_2}{2\pi} \; \cos \psi_1 \sin^3 \psi_1 = 0 ,
\end{equation}
and
\begin{equation} \nonumber
\int_{0}^{2\pi} \frac{{\rm d}\psi_1}{2\pi} \int_{0}^{2\pi} \frac{{\rm d}\psi_2}{2\pi} \; \sin i \psi_2 = 
\int_{0}^{2\pi} \frac{{\rm d}\psi_1}{2\pi} \int_{0}^{2\pi} \frac{{\rm d}\psi_2}{2\pi} \; \cos i \psi_2 = 0 ,
\end{equation}
we obtain
\begin{subequations} \label{eq:5.14}
\begin{align}
I_{1,c} & = - 2\gamma  \int_{0}^{2\pi} \frac{{\rm d}\psi_1}{2\pi} \; \cos \psi_1 \, \Phi(C_0\om_{L} \cos \psi_1 ) ,
\label{eq:5.14a} \\
I_{1,s} & = - \mu C_0 \int_{0}^{2\pi} \frac{{\rm d}\psi_1}{2\pi} \; \sin^2  \psi_1 +
\frac{A_3}{3} C_0^3 \int_{0}^{2\pi} \frac{{\rm d}\psi_1}{2\pi} \; \sin^4  \psi_1 
= - \frac{1}{2} \mu C_ 0 + \frac{A_3}{8} C_{0}^{3}.
\label{eq:5.14b}
\end{align}
\end{subequations}

\begin{figure}[!ht]
\centering 
\subfigure{\includegraphics[width=2.4in]{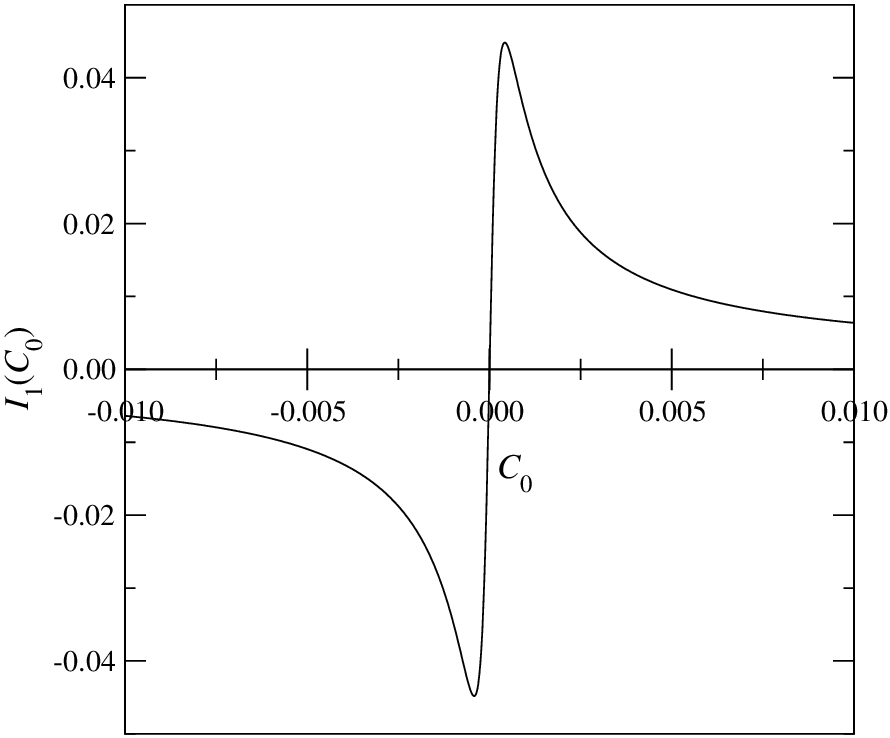}}
\hspace{2.0cm}
\subfigure{\includegraphics[width=2.4in]{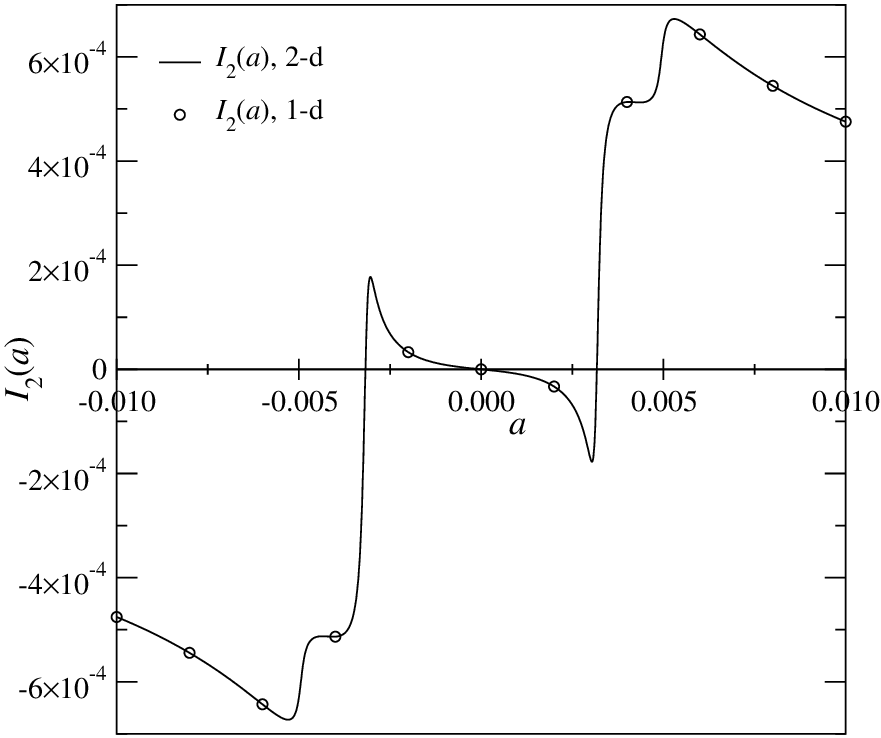}}
\caption{Left: plot of $I_1(C_0)$, defined in \eqref{eq:5.15}; right:
plot of $I_2(a)$, with the solid line corresponding to the $2$-dimensional integral defined in \eqref{eq:5.25a}
and the circles corresponding to the one-dimensional integral defined as the time average, in the sense of
\eqref{eq:5.22}, of the function $\cos \om t\,\Phi(a\cos\om t+b(nt))$.}
\label{figure:7}
\end{figure}

The two equations \eqref{eq:5.14} are coupled, since, even though $\mu$ does not appear explicitly in
\eqref{eq:5.14a}, the frequency $\om_{L}$ must satisfy the relation $\om^2=\om_{L}^2+\mu\,\epsilon$.
If $C_0\neq 0$, \eqref{eq:5.14b} gives $\mu=A_3 C_0^2/4$, while $C_0=0$ would not fix $\mu$ to any value.
We fix $C_0$ in such a way to make the integral
\begin{equation} \label{eq:5.15}
I_{1}(C_0) := \int_{0}^{2\pi} \frac{{\rm d}\psi_1}{2\pi} \; \cos \psi_1 \, \Phi(C_0\om_L \cos \psi_1 )
\end{equation}
vanish. One has
\begin{equation} \nonumber
I_{1}(0) = \Phi(0) \int_{0}^{2\pi} \frac{{\rm d}\psi_1}{2\pi} \; \cos \psi_1 = 0 ,
\end{equation}
which shows that $C_0=0$ is a zero of $I_1(C_0)$. To study the existence of other zeroes,
we compute numerically the integral $I_1(C_0)$, by approximating $\om_L$ with $\om$.
The result, given in Figure \ref{figure:7}, shows that there are no other zeroes.

This fixes the value of $C_0=0$ in the zeroth order approximation \eqref{eq:5.9}.
With $C_0=0$, the value of $\mu$ in \eqref{eq:5.14b} is left undetermined and can be set equal to $0$.
With the two values $C_0$ and $\mu$ being fixed, the right hand side of \eqref{eq:5.11}
becomes determined and $\Xi_1(t)$ in \eqref{eq:5.12} can be computed explicitly. One has
\begin{eqnarray} 
\Xi_{1}(t) & \!\!\! = \!\!\! & - \frac{2\gamma\e}{\om_{L}} 
\left( \sin \om_{L} t  \!\! \int_{0}^{t} \!\! {\rm d}  \tau \, \cos \om_{L} \tau \; \Phi(0) -
\cos \om_{L} t  \!\! \int_{0}^{t} \!\! {\rm d} \tau \, \sin \om_{L} \tau \; \Phi(0) \right) \nonumber \\
& \!\!\! + \!\!\! & \frac{2\e}{\om_{L}}  \sum_{i\in\II_{0}} \! B_{i} \!
\left( \sin \om_{L} t  \!\! \int_{0}^{t} \!\! {\rm d}  \tau \, \cos \om_{L} \tau \, \sin (in\tau) -
\cos \om_{L} t  \!\! \int_{0}^{t} \!\! {\rm d} \tau \, \sin \om_{L} \tau \, \sin (in\tau) \right) , \nonumber
\end{eqnarray}
so that, by computing the integrals explicitly, one finds
\begin{equation} \label{eq:5.16}
\Xi_{1}(t) = \frac{2\e\,B}{\om_{L}} \sin \om_{L} t + \beta (nt) + \rho \left( \cos \om_{L} t - 1 \right) ,
\end{equation}
where
\begin{equation} \label{eq:5.17}
B := \sum_{i\in\II_{0}} B_{i} \frac{in}{(in)^2-\om_{L}^2} , \qquad \rho := \frac{2\gamma\e}{\om_{L}^2} \Phi(0) ,
\qquad \beta(\psi) := - 2 \e  \sum_{i\in\II_{0}} \frac{B_{i}}{(in)^2-\om_{L}^2} \sin(i\psi) 
\end{equation}
and $\om_{L}$ has to be fixed at the value $\om$. Since
one has $\rho = 5.242 \times 10^{-5}$ and $2\e B/\om = -8.819 \times 10^{-3}$, for $\om_{L}=\om$, the last
contribution in \eqref{eq:5.16} is negligible, so that we can approximate
\begin{equation} \label{eq:5.18}
\xi_{1}(t) = X_{1}(\oo_L t) = \al \sin \om_{L} t + \beta(n t) , \qquad
\al:= C_1 + \frac{2\e B}{\om_{L}} ,
\end{equation}
In \eqref{eq:5.18} the constant $C_1$ is still arbitrary. In order to obtain the full expression for the first approximation,
we have to study the equation for the second approximation and impose the
requirement that the solution be bounded.

Of course, if we are interested in the first order contribution, we have to set $\om_{L}=\om$. However,
as explained before, if we want to compute the second order approximation we have to leave $\om_{L}$
as a free parameter, to be fixed together with $C_1$ by requiring the second order approximation to remain bounded.

\subsection{Second approximation} \label{sec:5.3}

The second approximation is obtained as the solution to the equation
\begin{equation} \nonumber 
\ddot  \xi = - \omega_{L}^2 \xi - \mu \e \xi_{1}(t) -
2 \e \sum_{i\in\II_{0}} B_{i} \sin( \xi_{1}(t) - in t) - 2 \gamma \, \e \Phi(\dot \xi_{1}(t)) + \frac{A_3 \e}{3} \, \xi^3_{1}(t) ,
\end{equation}
where the function $\xi_{1}(t)$ will be approximated by \eqref{eq:5.18} in the following. Therefore $\xi_2(t)$ is given by
\eqref{eq:5.11a} for $k=2$, where we write once more
\begin{equation} \nonumber
\bar\xi_2 \, \cos \om_{L} t + \bar y_2 \sin \om_{L} t = C_{2} \sin (\om_{L} t + \varphi_2) ,
\end{equation}
with $\varphi_2$ fixed to be zero and $C_2>0$ to be determined by imposing that no
secular terms appear when computing the third approximation. Therefore, we can write
the second approximation as
\begin{equation} \label{eq:5.20}
\xi_{2}(t) = C_2 \sin \om_{L} t + \Xi_{2}(t) , \qquad
\Xi_{2}(t) := \frac{\e}{\om_{L}} \Bigl( \sin \om_{L} t \,  M_{2,c}(\tau) - \cos \om_{L} t \,  M_{2,s}(\tau) \Bigr) .
\end{equation}
where $M_{2,c}(t)$ and $M_{2,s}(t)$ are defined in \eqref{eq:5.11b} with $k=2$.
For the function \eqref{eq:5.20} to be bounded one needs
\begin{equation} \label{eq:5.22}
I_{2,c} := \lim_{T\to+\io} \frac{M_{2,c} (T)}{T} = 0 , \qquad I_{2,s} := \lim_{T\to+\io} \frac{M_{2,s} (T)}{T} = 0 .
\end{equation}
We can write $I_{2,c} = I_{2,c,1} + I_{2,c,2} + I_{2,c,3} + I_{2,c,4}$ and
$I_{2,s} = I_{2,s,1} + I_{2,s,2} + I_{2,s,3} + I_{2,s,4}$, where

%
\begin{subequations} \label{eq:5.23}
\begin{align}
I_{2,c,1} & = 2 \sum_{i\in\II_{0}} \! B_{i} 
\int_{0}^{2\pi} \frac{{\rm d}\psi_1}{2\pi} \int_{0}^{2\pi} \frac{{\rm d}\psi_2}{2\pi} \; \cos \psi_1 \;
\sin (i\psi_2 - \al \sin \psi_1 - \beta(\psi_2) ) ,
\label{eq:5.23a} \\
I_{2,c,2} & = 
- 2 \gamma \int_{0}^{2\pi} \frac{{\rm d}\psi_1}{2\pi} \int_{0}^{2\pi} \frac{{\rm d}\psi_2}{2\pi} \; \cos \psi_1 \; 
\Phi(a \, \cos \psi_1 + b(\psi_2) ) ,
\label{eq:5.23b} \\
I_{2,c,3} & = 
- \mu \int_{0}^{2\pi} \frac{{\rm d}\psi_1}{2\pi} \int_{0}^{2\pi} \frac{{\rm d}\psi_2}{2\pi} \; \cos \psi_1 \; 
\left( \al \sin\psi_1 + \beta(\psi_2) \right) ,
\label{eq:5.23c} \\
I_{2,c,4} & = 
\frac{A_3}{3} \int_{0}^{2\pi} \frac{{\rm d}\psi_1}{2\pi} \int_{0}^{2\pi} \frac{{\rm d}\psi_2}{2\pi} \; \cos \psi_1 \; 
\left( \al \sin\psi_1 + \beta(\psi_2) \right)^3 ,
\label{eq:5.23d}
\end{align}
\end{subequations}
and, analogously,
\begin{subequations} \label{eq:5.24}
\begin{align}
I_{2,s,1} & = 2 \sum_{i\in\II_{0}} \! B_{i} 
\int_{0}^{2\pi} \frac{{\rm d}\psi_1}{2\pi} \int_{0}^{2\pi} \frac{{\rm d}\psi_2}{2\pi} \; \sin \psi_1 \;
\sin (i\psi_2 - \al \sin \psi_1 - \beta(\psi_2) ) ,
\label{eq:5.24a} \\
I_{2,s,2} & = 
- 2 \gamma \int_{0}^{2\pi} \frac{{\rm d}\psi_1}{2\pi} \int_{0}^{2\pi} \frac{{\rm d}\psi_2}{2\pi} \; \sin \psi_1 \; 
\Phi(a \, \cos \psi_1 + b(\psi_2) ) ,
\label{eq:5.24b} \\
I_{2,s,3} & = 
- \mu \int_{0}^{2\pi} \frac{{\rm d}\psi_1}{2\pi} \int_{0}^{2\pi} \frac{{\rm d}\psi_2}{2\pi} \; \sin \psi_1 \; 
\left( \al \sin\psi_1 + \beta(\psi_2) \right) ,
\label{eq:5.24c} \\
I_{2,s,4} & = 
\frac{A_3}{3} \int_{0}^{2\pi} \frac{{\rm d}\psi_1}{2\pi} \int_{0}^{2\pi} \frac{{\rm d}\psi_2}{2\pi} \; \sin \psi_1 \; 
\left( \al \sin\psi_1 + \beta(\psi_2) \right)^3 ,
\label{eq:5.24d}
\end{align}
\end{subequations}
with
\begin{equation} \nonumber
a := \al \om_{L} = C_1 \om_{L} + 2 B \e , \qquad 
b(\psi) := - 2 \e  \sum_{i\in\II_{0}} \frac{i n \, B_{i}}{(in)^2-\om_{L}^2} \cos(i \psi) ,
\end{equation}
%
%
%
Note that $b(nt)=\dot\beta(nt)$, where the derivative is with respect to $t$. In particular $\beta(\psi)$
is odd and hence $b(\psi)$ is even, so that $I_{2,s,2}=0$. Moreover $I_{2,c,3}=I_{2,c,4}=0$ and, 
since for $\al\neq0$ one has
\begin{equation} \nonumber
\cos \psi_1 \; \sin (i\psi_2 - \al \sin \psi_1 - \beta(\psi_2) ) = \frac{1}{\al}
\frac{{\rm d}}{{\rm d}\psi_1} \cos (i\psi_2 - \al \sin \psi_1 - \beta(\psi_2) ) ,
\end{equation}
also $I_{2,c,1}=0$. As a consequence one has $I_{2,c}=I_{2,c,2}$ and $I_{2,s}=I_{2,s,1}+I_{s,2,3}+I_{2,s,4}$.
The integrals $I_{2,s,3}$ and $I_{2,s,4}$ can be easily computed, up to corrections, and give
\begin{equation} \nonumber
I_{2,s,3} = -\mu \frac{\al}{2} , \qquad
I_{2,s,4} \approx \frac{A_3}{3} \times \frac{3\al^3}{8} = A_3 \frac{\al^3}{8} , 
\end{equation}
where terms of order $\al\|\beta\|_{\io}$ have been neglected in the latter, since $\|\beta\|_{\io}=1.644 \times 10^{-4}$.

In conclusion \eqref{eq:5.22} leads to the equations
\begin{subequations} \label{eq:5.25}
\begin{align}
I_{2}(a) := \int_{0}^{2\pi} \frac{{\rm d}\psi_1}{2\pi} \int_{0}^{2\pi} \frac{{\rm d}\psi_2}{2\pi} \; \cos \psi_1 \; 
\Phi(a \, \cos \psi_1 + b(\psi_2) ) & = 0 ,
\label{eq:5.25a} \\
2 \sum_{i\in\II_{0}} \! B_{i} 
\int_{0}^{2\pi} \frac{{\rm d}\psi_1}{2\pi} \int_{0}^{2\pi} \frac{{\rm d}\psi_2}{2\pi} \; \sin \psi_1 \;
\sin (i\psi_2 - \al \sin \psi_1 - \beta(\psi_2) ) + \frac{A_3\al^3}{8} - \frac{\mu \al}{2} & =  0 .
\label{eq:5.25b}
\end{align}
\end{subequations}

The integral $I_{2}(a)$ can be computed numerically, by approximating $\om_{L}=$ with $\om$.
The outcome is given in Figure \ref{figure:7}.
The results suggest the existence of three simple zeroes; besides $a=0$, there are two other zeroes
$a=\pm 3.174 \times 10^{-3}$. Since we are taking $C_1>0$ and hence $a>2\e B=-0.3116 \times 10^{-3}$,
only the positive zero has to be considered. The corresponding value of $C_1$, computed once more at $\om_{L}=\om$,  is
\begin{equation} \label{eq:5.26}
C_1 = \frac{a-2\e B}{\om} = \frac{3.174+3.116}{0.3534} \times 10^{-3} = 1.780 \times 10^{-2} ,
\end{equation}
which gives $\al=8.981 \times 10^{-3}$ in \eqref{eq:5.18}.

\begin{rmk} \label{rmk:5}
\emph{
In principle one should consider also the zero $a=0$ (which would give $C_1=8.817 \times 10^{-2}$).
However such a value makes $\mu$ disappear from the equations \eqref{eq:5.25} and hence does not fix $\om_L$.
To study the fate of such a solution, in particular to see whether it does correspond to a solution of the full equation,
one should go to higher orders, where we expect such a solution to disappear --- see also Section \ref{sec:7}.
Here we focus on the solution corresponding to $a\neq0$, since the latter correctly describes
the attractor relevant for the dynamics, as shown by the comparison below with the numerics.
}
\end{rmk}

\begin{figure}[H]
\centering 
\vspace{-.2cm}
\includegraphics[width=2.8in]{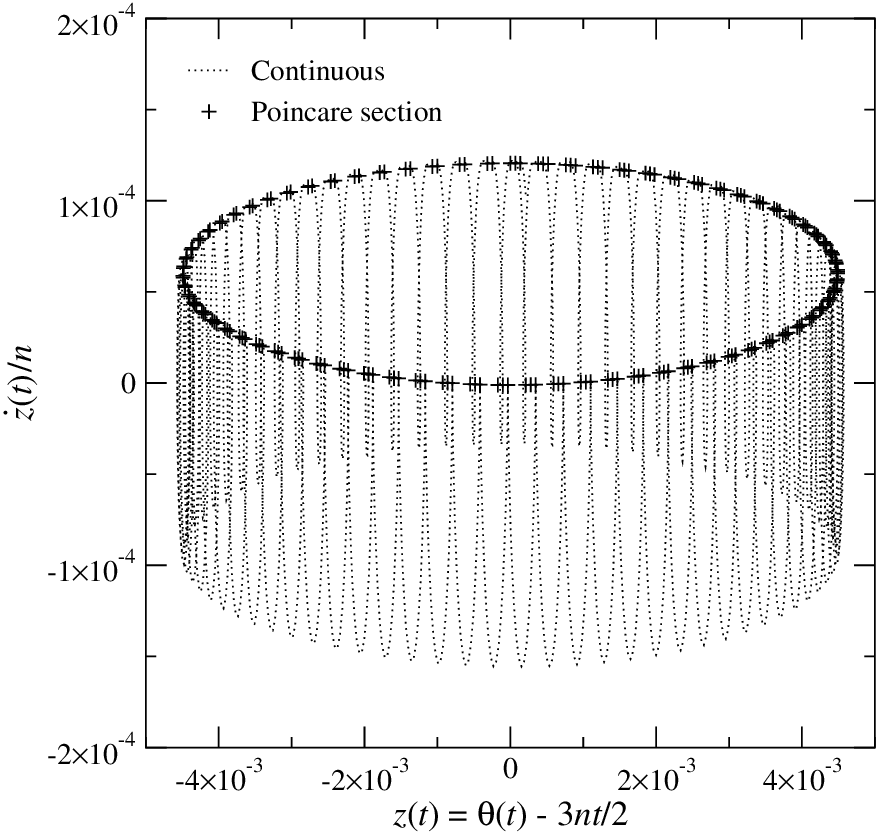}
\caption{Plot of $(z(t),\dot z(t)/n)$, with $C_1$ as in \eqref{eq:5.26} and $\om_{L}=\om$.
The crosses correspond to the Poincar\'e section obtained for $t$ an integer multiple of $T_0$, with $T_0=2\pi/n$.}
\label{figure:8}
\end{figure}

In \eqref{eq:5.25b} we can approximate the integral by expanding
\begin{eqnarray} 
& & \sin \left( i \psi_2 - \al \sin \psi_1 - \beta(\psi_2) \right) = \sin (i\psi_2) - \cos(i\psi_2) \left( \al \sin \psi_1 + \beta(\psi_2) \right)
\nonumber \\
& & \hspace{1cm} - \frac{1}{2} \sin (i\psi_2) 
\left( \al^2 \sin^2 \psi_1 + \beta^2(\psi_2) + 2 \al \sin \psi_1 \, \beta(\psi_2) \right) + \frac{1}{3!} \cos(i\psi_2) \left( \al \sin \psi_1 \right)^3 , 
\nonumber
\end{eqnarray}
where the other terms of order equal to or higher than three have been neglected.
The only non-zero contribution to the integral is
\begin{equation} \nonumber
J_2(a) := 4 \al \e \sum_{i,j\in\II_{0}} B_{i}
\int_{0}^{2\pi} \frac{{\rm d}\psi_1}{2\pi}  \int_{0}^{2\pi} \frac{{\rm d}\psi_2}{2\pi} \; \sin^2 \psi_1 \, \sin(i\psi_2) \,
\frac{B_{j}}{(jn)^2-\om_{L}^2} \sin(j \psi_2) ,
\end{equation}
where the explicit form of $\beta(\psi_2)$ has been used. Only the contributions with $j=\pm i$ are non-zero, so that
eventually we obtain
\begin{equation} \label{eq:5.27}
\frac{\mu\al}{2} - \frac{A_3 \al^3}{8} =
4 \al \e \sum_{i\in\II_{0}} \frac{B_{i}\left(B_{i}-B_{-i}\right)}{(in)^2-\om_{L}^2} 
\int_{0}^{2\pi} \frac{{\rm d}\psi_1}{2\pi}  \int_{0}^{2\pi} \frac{{\rm d}\psi_2}{2\pi} \; \sin^2 \psi_1 \, \sin^2(i\psi_2) =
\al \e \sum_{i\in\II_{0}} \frac{B_{i}\left(B_{i}-B_{-i}\right)}{(in)^2-\om_{L}^2} .
\end{equation}
An explicit computation, with $\om_{L}=\om$, gives
\begin{equation} \label{eq:5.28}
D:= \sum_{i\in\II_{0}} \frac{B_{i}\left(B_{i}-B_{-i}\right)}{(in)^2-\om^2} = 4.988 \times 10^{-4} .
\end{equation}
Inserting \eqref{eq:5.28} into \eqref{eq:5.27} yields
\begin{equation} \nonumber
\mu = \frac{A_3}{4} \al^2 + 2\e D = 1.319 \times 10^{-5} +
9.522 \times 10^{-5} = 1.084 \times 10^{-4} 
\end{equation}
and hence $\mu\e= 1.034 \times 10^{-5}$.
This fixes $\om_{L}$ to a value such that $\om^2-\om_{L}^2=\mu\e$,
so that $\om_L-\om \approx 1.46 \times 10^{-5}$; hence $\om_L$ is very close to $\om$.
For the value of $C_1$ in \eqref{eq:5.26} and $\om_L$ approximated by $\om$,
we consider the solution
\begin{equation} \nonumber
\theta(t) = \frac{3}{2} nt + \frac{1}{2} \xi_{1}(t) , \qquad
\frac{1}{2} \xi_1(t) = \frac{\al}{2} \sin \om  t - \e  \sum_{k\in\KK_{0}} \frac{A_{k}}{(k-3)^2 n^2 - \om^2} \cos((k-3)nt) ,
\end{equation}
with $\al/2=4.490 \times 10^{-3}$,
and define the libration as $z(t):=\theta(t)-(3/2)nt$. Then $(z(t),\dot z(t)/n)$ is as plotted in Figure \ref{figure:8}.

\subsection{Comparison with the numerical results} \label{sec:5.4}
\setcounter{equation}{0}

\begin{figure}[!ht]
\centering 
\subfigure{\includegraphics[width=2.6in]{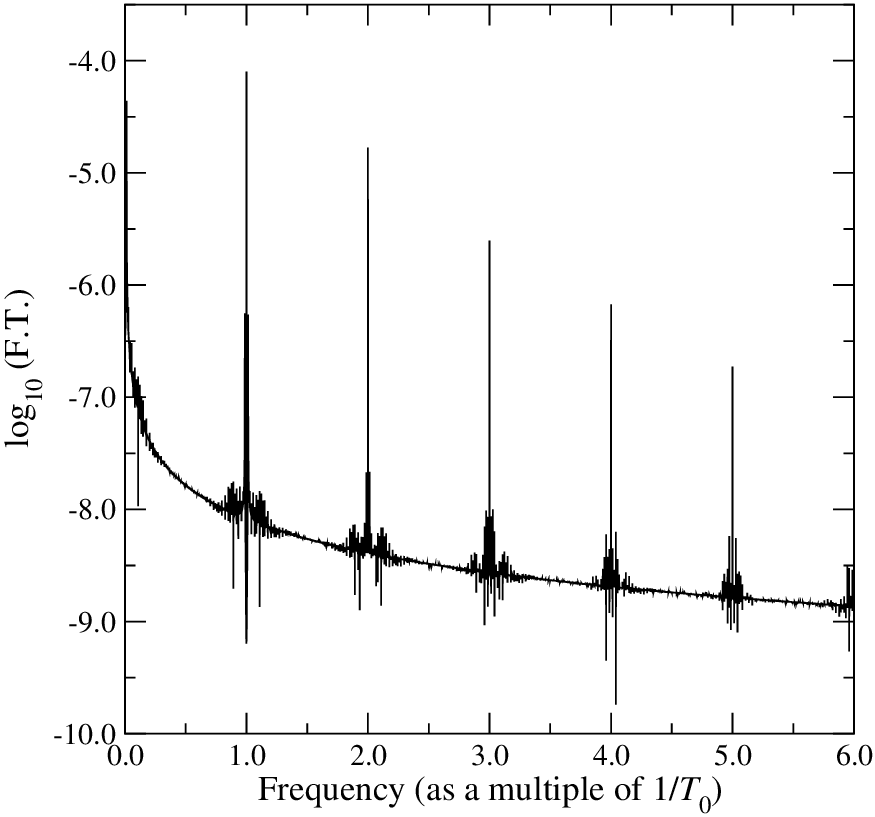}}
\hskip1cm
\subfigure{\includegraphics[width=2.6in]{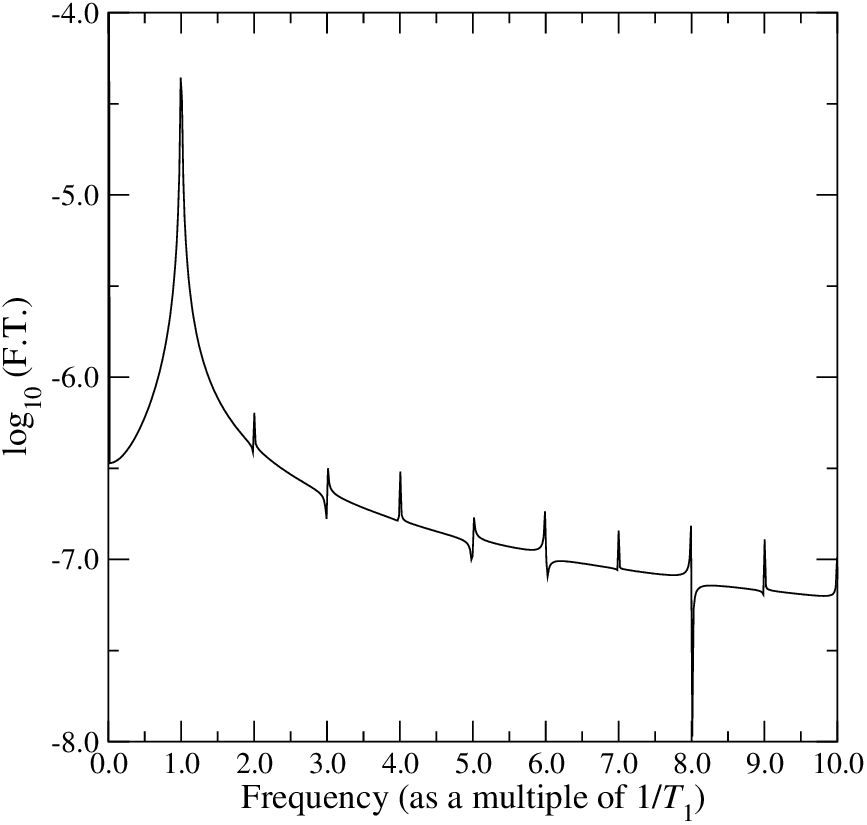}}
\caption{Fourier transform of $\dot\theta$ for the numerical quasi-periodic solution of Figure \ref{figure:2};
here $T_0=2\pi/n$ and $T_1=73.9034 \, T_0=2\pi/\om_L$. There are peaks at integer multiples of
$1/T_0$ and $1/T_1$, respectively.}
\label{figure:9}
\end{figure}

The approximate solution found in the previous sections has to be compared with that found by numerical analysis,
see Figure \ref{figure:2}. Fast Fourier Transform analysis produces the spectrum of the derivative of the
numerical solution in Figure \ref{figure:9}: the figure to the left, on a larger scale, shows the presence 
of peaks at multiples of the frequency $n$,
while the figure to the right, on a smaller scale, shows that there are peaks as well at multiples of the
frequency $\om_L=n/73.9034 = 0.3530$. Therefore, with respect to the approximate analytical solution,
one has $\om - \om_L = 4 \times 10^{-4}$ and hence
$\om^2-\om_{L}^2=0.3534^2-0.3530^2=2.8256 \times 10^{-4}$, which is larger than the analytical value.
However, apart from that, the agreement between analytical approximations and numerical results is very good.
Indeed, we expect such a difference to be negligible as far as we are interested in the amplitude of the solution.
Nevertheless the phase shift will become appreciable on longer timescales.

It seems likely to us that, in order to obtain the correct value of the frequency $\om_L$,
one has to go to higher orders of approximation: because of the rapid variations of the function $\Phi$,
the computations depend very sensitively on the approximation $\dot\xi_{k-1}$
appearing in its argument in \eqref{eq:5.11a}. More precisely, when studying the third approximation
$\xi_{3} = C_3 \sin \om_L t + \Xi_{3}(t)$, we expect the following to happen:
requiring the solution to remain bounded should fix $C_2$ to a value very close to $C_1$,
while the corrections $\mu$ to the slow frequency should be appreciably larger than the value
found at the previous step. Then the second approximation $\xi_2(t)$ in \eqref{eq:5.20}
should provide not only the right amplitude, but the right frequency as well --- see Section \ref{sec:7}
for further comments.

\begin{rmk} \label{rmk:6}
\emph{
It is important to point out that, even though we are using the term `quasi-periodic' for the attractor
that we have studied in this section,
in principle the two frequencies $\om_L$ and $n$ could be commensurate. In fact, since we are determining
$\om_L$ numerically, it is not possible to exclude the possibility that the ratio $\om_L/n$ might be a rational number.
A better definition could be `multi-periodic' solution, since the solution appears as a superposition
of two oscillations, each with its own frequency. However, if we let the parameters change slightly,
we find the ratio to be irrational most of the time
(see Section \ref{sec:6} below). Thus, we can conclude that the attractor is very likely
to be genuinely quasi-periodic.
Note also that, even if the ratio were rational and hence the solution were periodic,
its period would be very large, since in any case the ratio would not be close to any rational $p/q$, with $q$ small.
}
\end{rmk}

\section{Numerical study of the stability of the attractors} \label{sec:6}
\setcounter{equation}{0}

The stability of the periodic attractors for the system described by \eqref{eq:2.1} can be studied by considering the
corresponding Poincar\'e map and computing the eigenvalues of the linearised system around its fixed points~\cite{MC,GH}.

As shown in Section \ref{sec:3}, there are four periodic orbits for each resonance and hence four fixed points for
a suitable iteration of the Poincar\'e map. Among the fixed points, for $\e$ small enough, two are stable and two are unstable.
We can confine ourselves to the fixed points corresponding to the values of the initial phase
of the first approximation $\Th_1^{(1)}$ and $\Th_1^{(2)}$; indeed, by symmetry, the points corresponding
to the phases $\Th_1^{(3)}$ and $\Th_1^{(4)}$ have the same stability of 
$\Th_1^{(1)}$ and $\Th_1^{(2)}$, respectively.

For the values of the parameters given in Section \ref{sec:2}, we find the eigenvalues given in Table \ref{table:4}.
Numerically, one observes that, when the periodic attractors are both unstable, quasi-periodic attractors such as that studied
in Section \ref{sec:5}, appear. This happens for the resonances 1:2, 3:2 and 2:1, besides the retrograde resonances
which, however, do not play a relevant role in the dynamics of the system. The stable
quasi-periodic attractors appear by a Hopf bifurcation, occurring when changing the parameters
of the system \cite{MC}.

As in Section \ref{sec:5}, in the following we concentrate on the resonance 3:2, which is the most interesting one
in the case of Mercury, but a similar analysis could be easily extended to the other cases.
However, as a comparison with the results in Table \ref{table:4} shows, only quasi-periodic attractors
with $\dot\theta/n \approx 1/2$, $3/2$ and $2$  are found to exist among the prograde ones,
so that we could confine ourselves to these three cases.

\begin{table}[h]
\centering
\vspace{.2cm}
\begin{tabular}{|l|l|l|r|} \hline
\multicolumn{1}{|c}{Resonance} & \multicolumn{1}{|c}{$\theta(0)$} &
\multicolumn{1}{|c}{{\raisebox{2.7ex}{ }}{\raisebox{-1.1ex}{
}}$\dot\theta(0)/n$} &
\multicolumn{1}{|c|}{Eigenvalues} \\ \hline
$-1:1$ & 0.36289190131044645472  & $-1.00004242365858443089$ & $2.024\times 10^{-9}$\\
$-1:2$ & 0.14580421354300878946  & $-0.50006849412399051400$ & $2.783\times 10^{-9}$\\ \hline
1:2  & 3.14129563170348761883  &   0.49980635331803679181 & $0.9669, 1.0342$\\
1:1  & 3.14151499384565687042  & 0.99986201340697665762 & $-4.461\times 10^{-4}$\\
3:2  & 3.14150380436395113505  & 1.50005973350740330252 & $1.055\times 10^{-4}$\\
2:1  & $3.26027930307144126711\times 10^{-5}$ & 2.00012557558534916792 & $1.786\times 10^{-3}$\\
5:2  & 3.14140519201664595044  & 2.50012075040501328073 & $-3.628\times 10^{-4}$\\
3:1  & 3.14109199137670843320  & 3.00009814397107114853 & $-2.636\times 10^{-5}$\\
7:2  & 3.14027640889440704126  & 3.50007711111008245662 & $-3.835\times 10^{-6}$\\
4:1  & 3.13797190712320535390  & 4.00006157245270746253 & $-6.337\times 10^{-7}$\\ \hline\hline
$-1:1$ & 1.20792006104664609582  & $-0.99995757575039987029$ & $0.9992, 1.0008$\\
$-1:2$ & 1.42500286020411552411  & $-0.49993150584317966306$ & $0.9986, 1.0014$\\ \hline
1:2  & 1.57112385469851460569  &   0.50019364882055637631  & $2.342\times 10^{-6}$\\
1:1  & 1.57068938450889863242  &   1.00013792675908729505  & $0.9048, 1.1042$\\
3:2  & 1.57075984135159670901  &   1.49994030293249049891  & $0.9185, 1.0889$\\
2:1  & 1.57099968204819540739  &   1.99987444617026058657  & $0.9433, 1.0638$\\
5:2  & 1.57101812013673537458  &   2.49987925336853351100  & $0.9613, 1.0395$\\
3:1  & 1.57130265033260668261  &   2.99990185551468461907  & $0.9760, 1.0246$\\
7:2  & 1.57211353266178141100  &   3.49992288814339436128  & $0.9854, 1.0147$\\
4:1  & 1.57441706775605802984  &   3.99993842708145177608  & $0.9914, 1.0087$\\ \hline
\end{tabular}
\caption{The initial conditions for the existing periodic solutions and
the corresponding eigenvalues. The eigenvalues $\lambda_{1, 2}$ for each periodic solution
either form a conjugate pair, in which case $|\lambda_{1, 2}|-1$ is given, or are both real, in which case both are given.}
\label{table:4}
\end{table}

Writing~\eqref{eq:2.1} as
\begin{equation} \label{eq:6.1}
\ddot  \theta = -  \zeta \, A_{3} \sin (2\theta - 3nt) - S \, \zeta \sum_{k\in\KK_0} A_{k} \sin( 2 \theta - k n t) - 
\lambda \eta \, F(\dot \theta) ,
\end{equation}
where $S=\lambda=1$; in terms of $\xi$, defined according to \eqref{eq:5.1}, the equation becomes
\begin{equation} \label{eq:6.2}
\ddot  \xi = -  \zeta \, A_{3} \sin \xi - S \, \zeta \sum_{i\in\II_0} B_{i} \sin( \xi - i n t) - 
\lambda \eta \, \Phi(\dot \theta) ,
\end{equation}
where the set $\II_0$ and the coefficients $B_i$ are defined as after \eqref{eq:5.3}.
We study the transitions in the dynamics of the system described by \eqref{eq:6.1}
when we vary either the parameter $\lambda$ (at fixed $S=1$) or the parameter $S$ (at fixed $\lambda=1$).

\begin{figure}[H]
\centering 
\subfigure{\includegraphics[width=2.7in]{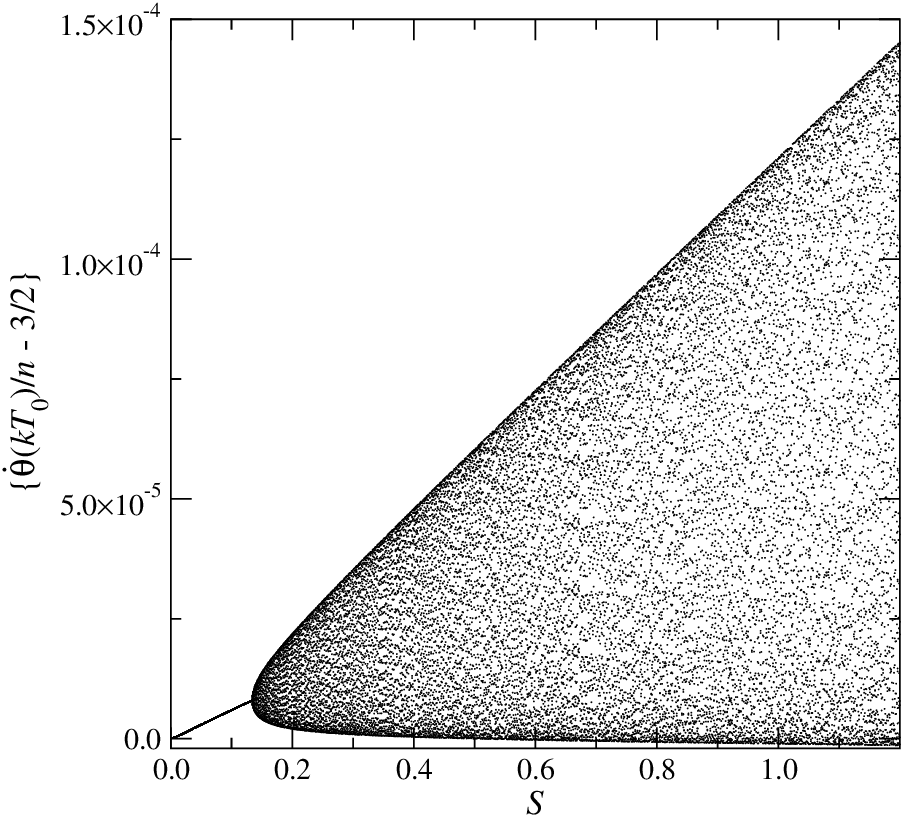}}
\hskip1cm
\subfigure{\includegraphics[width=2.8in]{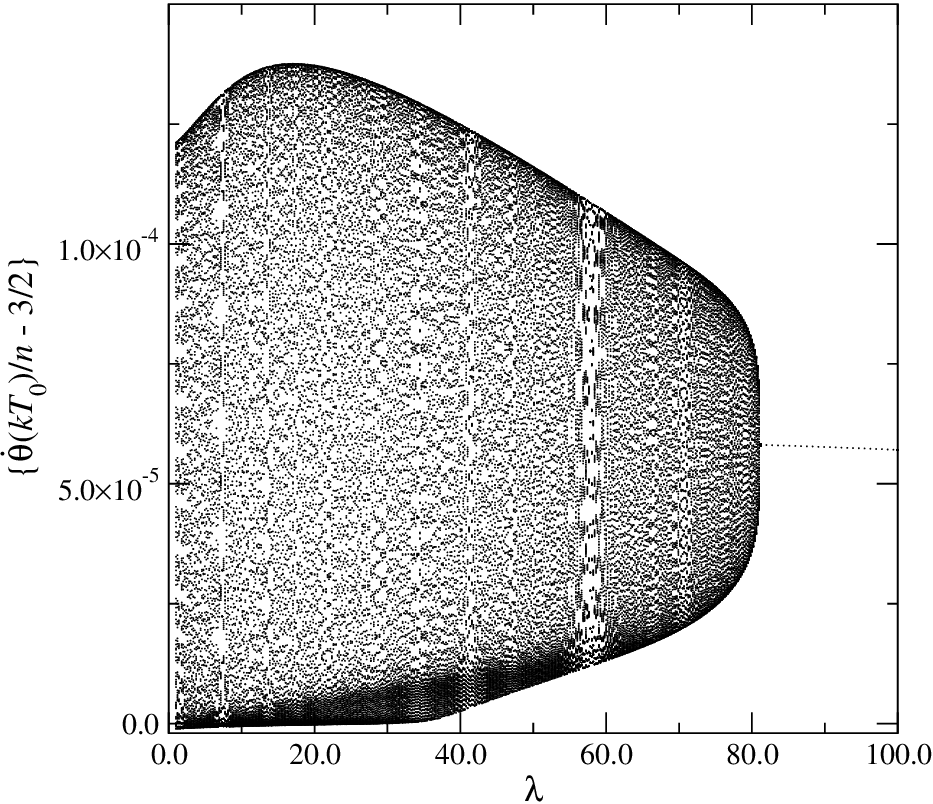}}
\caption{Bifurcation diagram with parameter $S$ (left) and $\lambda$ (right).
Both diagrams have been produced by plotting a set of values of $\dot\theta(kT_0)/n-3/2$, with $T_0=2\pi/n$,
for $k$ in a suitable large set of integers.}
\label{figure:10}
\end{figure}

The bifurcation diagram with the parameter $S$ in Figure \ref{figure:10} shows that the periodic solution 
with frequency $3n/2$ is stable up to the value
$S \approx 0.134$, where the solution loses stability and a stable quasi-periodic solution appears.
By increasing $S$, the amplitude of the oscillations increases as well. At $S=1$ we have
the quasi-periodic solution studied in Section \ref{sec:5}, while the periodic solution with frequency $3n/2$
is unstable. The numerical investigation in \cite{BDG2} demonstrates that the solution with velocity $\dot\theta$ close to  $3n/2$
is the main attractor for the values of the parameters as in Section \ref{sec:2}: more than $42\%$ of the initial conditions
$(\theta,\dot\theta) \in [0,2\pi] \times [0,5n]$ are captured by such an attractor.

\begin{figure}[H]
\centering 
\subfigure{\includegraphics[width=2.4in]{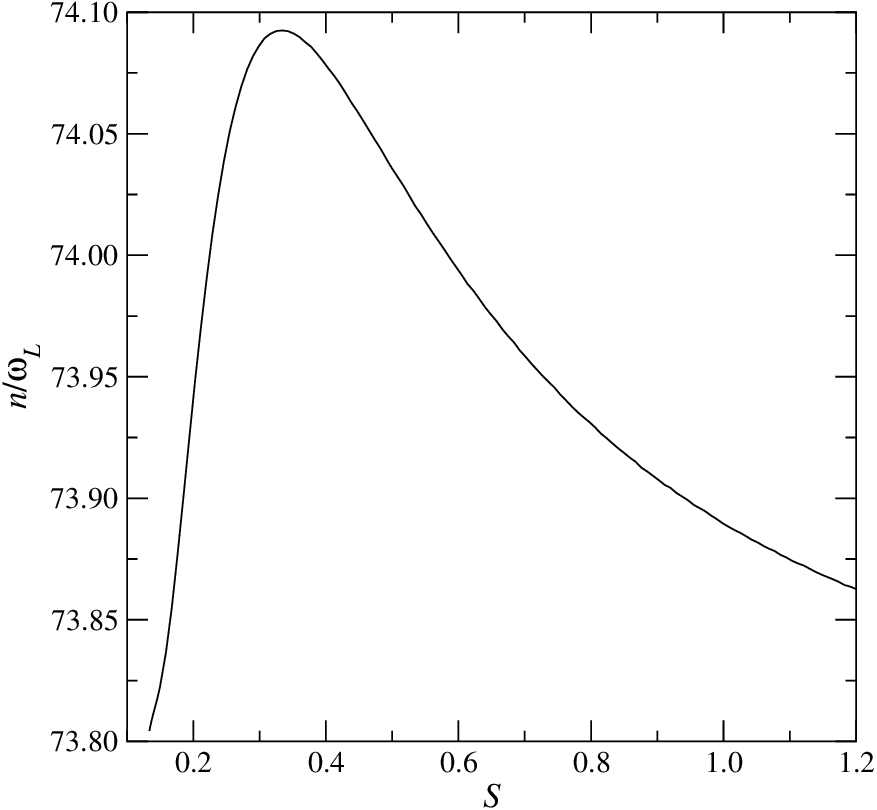}}
\hskip2cm
\subfigure{\includegraphics[width=2.35in]{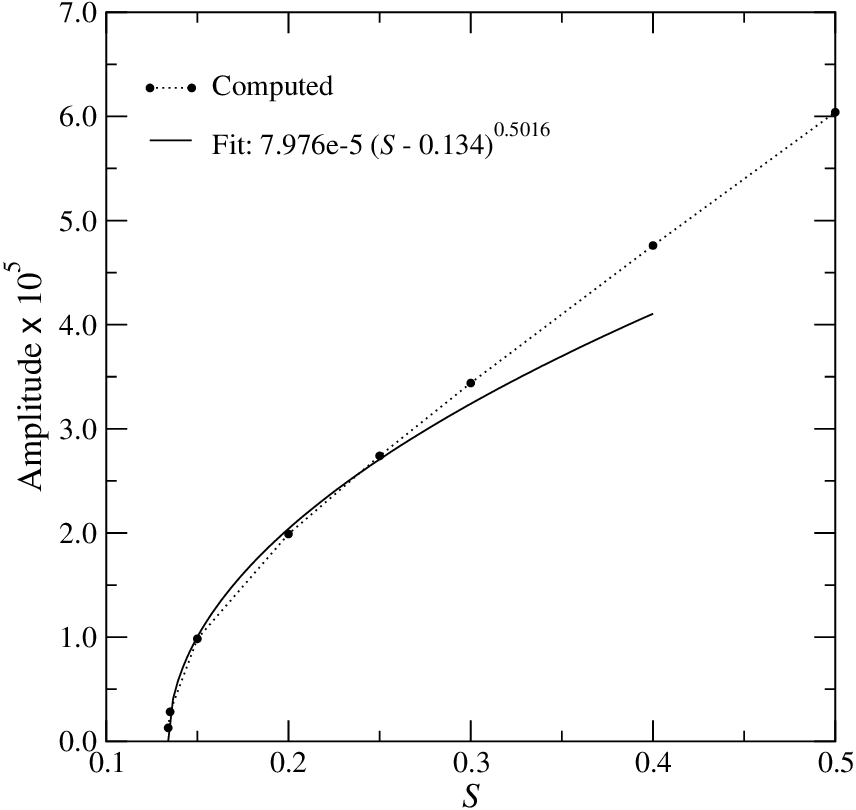}}
\caption{Low frequency $\om_{L}$ (left) and amplitude $A$ (right) of the quasi-periodic solution 
versus the bifurcation parameter $S$.}
\label{figure:11}
\end{figure}

In principle, the attractor studied in Section \ref{sec:5}, which we call quasi-periodic, could be periodic:
indeed, if $\om_L/n$ is a rational number $P/Q$,
then the trajectory closes after a suitable time $T$ (large, since $Q$ would be large).
Hence, we investigate
numerically how the slow frequency $\om_L$ changes when varying a parameter of the system.
In Figure \ref{figure:11} we plot the non-fixed frequency $\om_{L}$ of the
quasi-periodic solution as a function of the parameter $S$. The apparent
continuity of the curve suggests that,
up to a zero-measure set of values of $S$, the two frequencies $\om$ and $n$ are
incommensurate, so that the motion is genuinely quasi-periodic.

On the Poincar\'e section obtained by sampling the system evolution at integer multiples of $T_0=2\pi/n$,
the periodic solution appears as a finite set of points (3 for the 3:2 resonance). By contrast,
the quasi-periodic solution has support on a torus which intersects the section along a closed curve.
In terms of the dynamics on the section we have a Hopf bifurcation, with the transition from a fixed point to
a closed curve \cite{MC}. A plot of the amplitude $A$ of the quasi-periodic solution versus the parameter $S$
is also given in Figure \ref{figure:11}: a fit obtained by looking at the leftmost points in the figure gives
$A=A_0(S-S_0)^{\kappa}$, with $A_0=7.976 \cdot 10^{-5}$, $S_0=0.134$ and $\kappa=0.5016$; the numerical value of the
exponent $\kappa$ is in agreement with the critical exponent $1/2$ typical of the Hopf bifurcation \cite{MC}.

From a physical point of view, the bifurcation diagram with parameter $S$ in Figure \ref{figure:10} can be interpreted as follows.
For $S=0$ there is a fixed point at the origin, which is asymptotically stable because of the presence of dissipation.
When $S$ becomes large enough (that is, if the forcing terms not included in the unperturbed system can not
be neglected any longer), the fixed point loses stability and an asymptotically stable
closed curve appears nearby by Hopf bifurcation.

A similar phenomenology is observed at $S=1$, when varying the parameter $\lambda$; see Figure \ref{figure:10}.
Actually we need a dissipation much larger than the physical value $\lambda=1$. Again the transition from
the periodic to the quasi-periodic solutions is described by a Hopf bifurcation for the dynamics of the
corresponding Poincar\'e section.
The bifurcation diagram with the parameter $\lambda$ in Figure \ref{figure:10} may also be expected on physical grounds.
If there is no dissipation, the Poincar\'e map corresponding to the equation \eqref{eq:6.2} has both a fixed point
near the origin and a large measure of invariant tori encircling such a point \cite{AKN,MNT,BC}.
As soon as $\lambda>0$, all but one of the tori are destroyed, this one being attractive; the fixed point exists as well
but it is unstable. By taking larger values of $\lambda$, at some point the torus is destroyed, while the fixed point
becomes stable and attracts all trajectories starting from initial data nearby.

\section{Dynamics far from the attractors} \label{sec:4}
\setcounter{equation}{0}

We now look at the pre-capture dynamics, that is, the period during which
the satellite is decelerating but before it has been captured. This will
enable us to estimate times to capture, and we accomplish this by making
approximations that greatly simplify the dissipation term.

Starting from \eqref{eq:2.1}, and assuming that (i) $\dot \theta > n$ and (ii) $\dot\theta$ is not close to any kink,
we make the approximation $\eta \,F(\dot \theta) = a - b\dot \theta$. This yields $a,b > 0$ --- see \cite[Figure 3]{BDG2}.
For $\dot \theta < n$, $a$ would be negative, but the argument below, suitably adapted, would still work.
In practice, both $a$ and $b$ are small. For instance, expanding $F(\dot \theta)$ 
around $\dot \theta = 1.75 n$ gives $a \approx 1.1\times 10^{-5}$ and $b\approx 1.3\times 10^{-7}$.

It is convenient to rescale time by $\tau = b t$, so that $\dot \theta = b \theta'$, where the prime
denotes the derivative with respect to $\tau$.
Using the above, we can approximate \eqref{eq:2.1} as
\begin{equation} \label{eq:4.1}
\theta'' = -\frac{a}{b^2} + \theta' - \frac{\zeta}{b^2} \sum_{k\in\KK}  A_k \sin \left( 2\theta  - \frac{k n}{b}\tau\right) .
\end{equation}
%

\subsection{The splitting argument}

We now make an estimate of the rate at which $\dot \theta$ decreases with time, 
over long time scales, typically of the order of $10^6 T_0$, where $T_0 = 2\pi/n$.
We start by defining $s(\tau)$ to be the solution of $s'' = s' -a/b^2 $
with initial conditions $s(0) = s_0$ and $s'(0) = s_1$, which gives
\begin{equation} \nonumber
s(\tau) = s_0 - s_1 +  \frac{a}{b^2} + \frac{a}{b^2}\tau +
e^\tau\left(s_1 -\frac{a}{b^2}\right) =  s_0 + s_1\tau + O\left(\tau^2\right).
\end{equation}

\begin{figure}[ht]
\centering
\subfigure{\includegraphics*[width=2.4in]{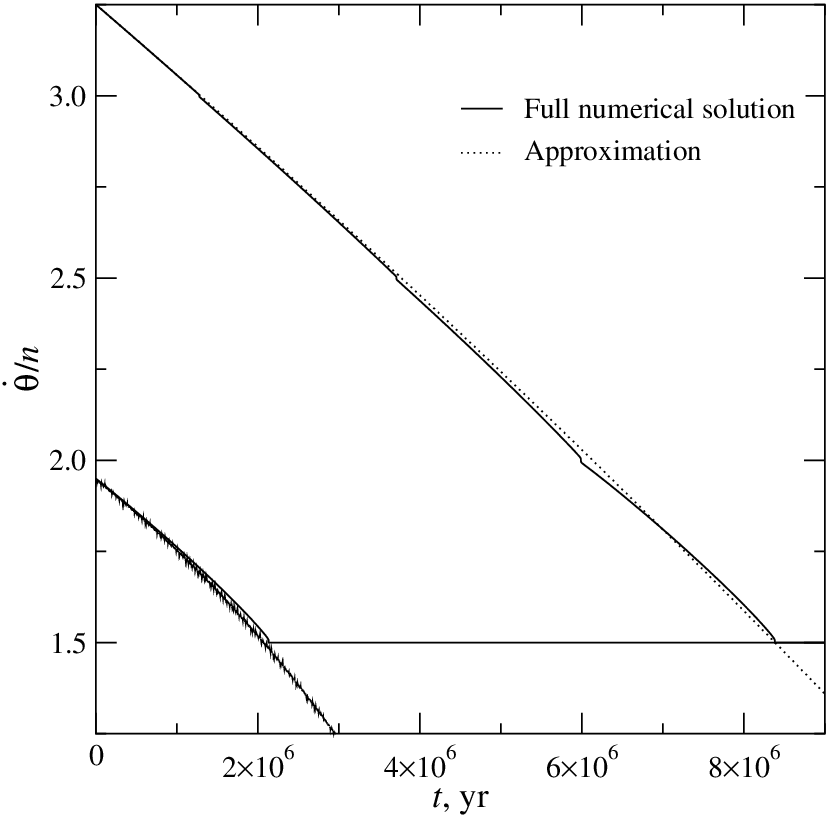}}
\hspace{2cm}
\subfigure{\includegraphics*[width=2.6in]{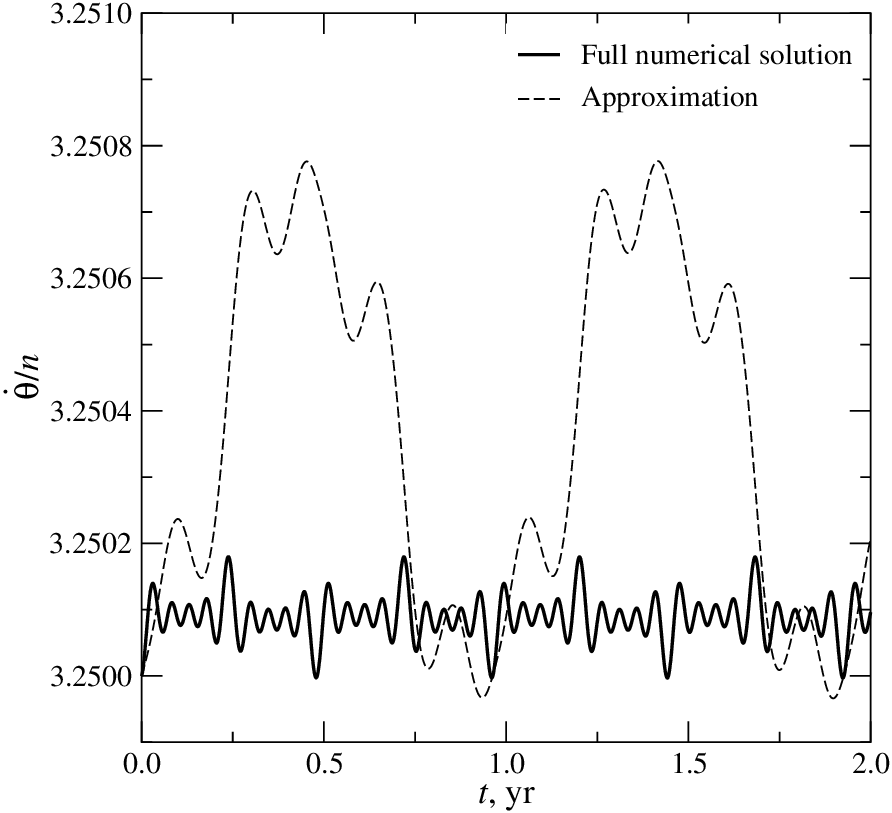}}
\caption{Left: The approximation in \eqref{eq:4.9} for  $\dot \theta(t)$ (dashed lines),
compared with the full numerical solution (solid lines), with initial spin rate $\do\theta(0)=3.25 n$ and
$\dot\theta(0)=1.95n$; right:
the same, but on a small time scale, with initial spin rate $\dot \theta(0) = 3.25 n$ and $\theta(0) = 1.7$.}
\label{figure:5}
\end{figure}

We write $\theta(\tau) = s(\tau) + f(\tau)$, in which $s(\tau)$ accounts for
the slow decay of $\dot \theta$ apparent on a large time scale and $f(\tau)$
represents a small amplitude, high frequency correction term,
which is visible only on a smaller time scale --- see Figure \ref{figure:5}. In other words,
we split $\theta(\tau)$ into fast and slow components. Substituting this expression
for $\theta(\tau)$ in \eqref{eq:4.1}, and bearing in mind the ODE obeyed by $s(\tau)$, we find that
\begin{equation} \nonumber 
f'' = f' - \frac{\zeta}{b^2}\sum_{k\in\KK} A_k\sin\left(\omega_{k} \tau + \phi + 2f(\tau)\right),
\end{equation}
where $\phi = 2 s_0$,
%
$\omega_{k} := 2s_1 - k n/b$
%
and where we have used the Taylor series to $O\left(\tau^2\right)$ for $s(\tau)$. 

We now make the following assumption: since $\omega_{k}$ is $O(1/b)$ and is therefore large,
we neglect the term $2f(\tau)$ and so obtain
\begin{equation} \label{eq:4.3}
f'' = f' - \frac{\zeta}{b^2}\sum_{k\in\KK} A_k\sin\left(\omega_{k} \tau + \phi \right) ,
\end{equation}
for an approximation to the ODE that defines $f(\tau)$. Conveniently, this
ODE is linear with constant coefficients, and so is straightforward to solve.
With initial conditions $f(0) = f_0$ and $f'(0) = f_1$, we find
\begin{multline} \label{eq:4.4}
f(\tau) = f_0 + f_1\left(e^\tau - 1\right) -
\frac{\zeta}{b^2} \sum_{k\in\KK} \frac{A_k}{1+\omega_{k}^2}\left[\frac{1}{\omega_{k}}\cos(\omega_{k}\tau + \phi)
- \sin(\omega_{k} \tau + \phi) \right.\\ \left.
+ \; e^\tau(\omega_{k} \cos\phi + \sin\phi) - (1+\omega_{k}^2)\frac{\cos\phi}{\omega_{k}}\right].
\end{multline}
We now consider the initial conditions, noting that we are free to choose initial
values for $\theta(0)=\theta_0 = s_0 + f_0$ and $\theta'(0)=\theta_1 = s_1 + f_1$. Once $\theta_0$, $\theta_1$
are specified, any values of $s_0$, $s_1$, $f_0$ and $f_1$
that satisfy these constraints can be chosen.

Recall that in order to derive \eqref{eq:4.3}, we assumed that $|f(\tau)| \ll 1$. To be consistent with this, we therefore
choose $f_0 = 0$, from which we immediately deduce that $s_0 = \theta_0$, so that $\phi = 2\theta_0$.

We need to take a little more care over the choice of $f_1$. Since
$\omega_{k}$ is $O(1/b)$, one has $|\omega_{k}|\gg 1$. Hence, in order to guarantee the
smallness of $f(\tau)$, we need if possible to choose $f_1$ so as to cancel out the
largest terms in the sum in \eqref{eq:4.4}, these being
$\omega_{k} e^\tau\cos 2 \theta_0$ and $-\omega_k\cos 2 \theta_0$, both of which are $O(1/b)$.
This cancellation can be accomplished by setting
\begin{equation} \nonumber
f_1 = \frac{\zeta}{b^2} \, \cos2 \theta_0 \sum_{k\in\KK}\frac{A_k \omega_k}{1+\omega_k^2},
\end{equation}
whereupon the expression for $f(\tau)$ becomes
\begin{equation} \label{eq:4.5}
f(\tau) = - \frac{\zeta}{b^2}\sum_{k\in\KK} \frac{A_k}{1+\omega_k^2}
\left[ \frac{\cos(\omega_k\tau + 2 \theta_0) - \cos 2 \theta_0}{\omega_k} -
\sin(\omega_k\tau + 2 \theta_0) + e^\tau \sin2 \theta_0\right].
\end{equation}
Since $\tau = b t$ with $b\sim 10^{-7}$ and the timescale that we consider 
is $t\sim 10^6$, $e^\tau$ is $O(1)$.

Finally we consider $s_1$. Since $\theta_1$ is given and we have already specified $f_1$, we immediately find that
\begin{equation} \nonumber
s_1 = \theta_1 - f_1 = \theta_1 - \frac{\zeta}{b^2} \, \cos 2 \theta_0 \sum_{k\in\KK}
\frac{A_k \omega_{k}}{1+\omega_{k}^2}.
\end{equation}

Hence, since $\theta(\tau) = s(\tau) + f(\tau)$, we have
\begin{multline} \label{eq:4.6}
\theta(\tau) = \theta_0 + \frac{a}{b^2}\tau + \frac{e^\tau-1}{b^2}\left[b^2 \theta_1 - a - \zeta
\cos 2 \theta_0\sum_{k\in\KK} \frac{A_k \omega_{k}}{1+\omega_{k}^2}\right] \\
- \frac{\zeta}{b^2}\sum_{k\in\KK} \frac{A_k}{1+\omega_{k}^2}\left[\frac{\cos(\omega_{k}\tau +
2 \theta_0) - \cos 2 \theta_0}{\omega_{k}} -
 \sin(\omega_{k}\tau + 2 \theta_0) + e^\tau \sin2 \theta_0\right].
\end{multline}

We now approximate this expression, starting from the fact that $| \omega_k| \gg 1$.
Additionally, we express everything in terms of $t$ rather than $\tau$, where $\tau = bt$, so 
that $b \theta_1 = \dot \theta(0)$. Then \eqref{eq:4.6} becomes
\begin{equation} \label{eq:4.7}
\theta(t) \approx \theta_0 + \frac{at}{b} + \frac{(e^{bt}-1)}{b}\left[
\dot \theta(0) - \frac{a}{b} - \zeta \cos 2 \theta_0 \sum_{k\in\KK} \frac{A_k}{\Omega_k}\right]
- \zeta \sum_{k\in\KK} \frac{A_k}{\Omega_k^2}\left[e^{bt} \sin 2 \theta_0 - \sin(\Omega_k t + 2 \theta_0)\right],
\end{equation}
where $\Omega_k = b\omega_k$, so
\begin{equation} \label{eq:4.8}
\Omega_k = \dot{s}_0 - k n = \dot \theta(0) -k n - \zeta \cos 2 \theta_0\sum_{k\in\KK} \frac{A_k}{\Omega_k}.
\end{equation}
Note that this equation defines $\Omega_k$ implicitly --- this
is a consequence of the way the initial conditions have been assigned.
For $n$  sufficiently large, the approximation $\Omega_k \approx \dot \theta(0) - kn$ will be good.
We investigate this approximation in practice in the next subsection.

Differentiating \eqref{eq:4.7} we find
\begin{equation} \label{eq:4.9}
\dot \theta(t) \approx \frac{a}{b} + e^{bt} \left[\dot \theta(0) -\frac{a}{b} \right] +
\zeta \sum_{k\in\KK}  \left[ \frac{A_k}{\Omega_k} \cos(\Omega_k t + 2 \theta_0) 
- e^{bt} \cos 2 \theta_0 \frac{A_k}{\Omega_k} \right] ,
\end{equation}
where we have neglected a term of order $b$ in the sum.

\subsection{The approximations in practice} \label{sec:4.2}

Several assumptions have been used to derive \eqref{eq:4.7} and \eqref{eq:4.9}, 
so we now investigate numerically how good these approximations are in practice.

First of all, we compare `exact' numerical solution $\Omega_k$ to \eqref{eq:4.8}
with the approximation $\Omega_k \approx \Omega_k^{\rm app} := \dot \theta(0) - k n$, for the
two values $\dot \theta(0) = 1.95 n$ and $3.25 n$. We find

\begin{itemize}
\itemsep0em
\item For $\dot \theta(0) = 1.95 n$, $\max_{\theta_0\in\mathbb T,\, k\in\QQ} |\Omega_k - \Omega_k^{\rm app}| \approx 0.073$;
\item For $\dot \theta(0) = 3.25 n$, $\max_{\theta_0\in\mathbb T,\, k\in\QQ} |\Omega_k - \Omega_k^{\rm app}| \approx 0.010$.
\end{itemize}

Note that $\min_{\theta_0\in\mathbb T,\,k\in\QQ} \left|\Omega_k\right| \approx 6.5$, so we are justified
in using the approximation \eqref{eq:4.8}, since it leads to a relative error of no more than 1.1\%.

\begin{figure}[!ht]
\centering
\includegraphics*[width=2.6in]{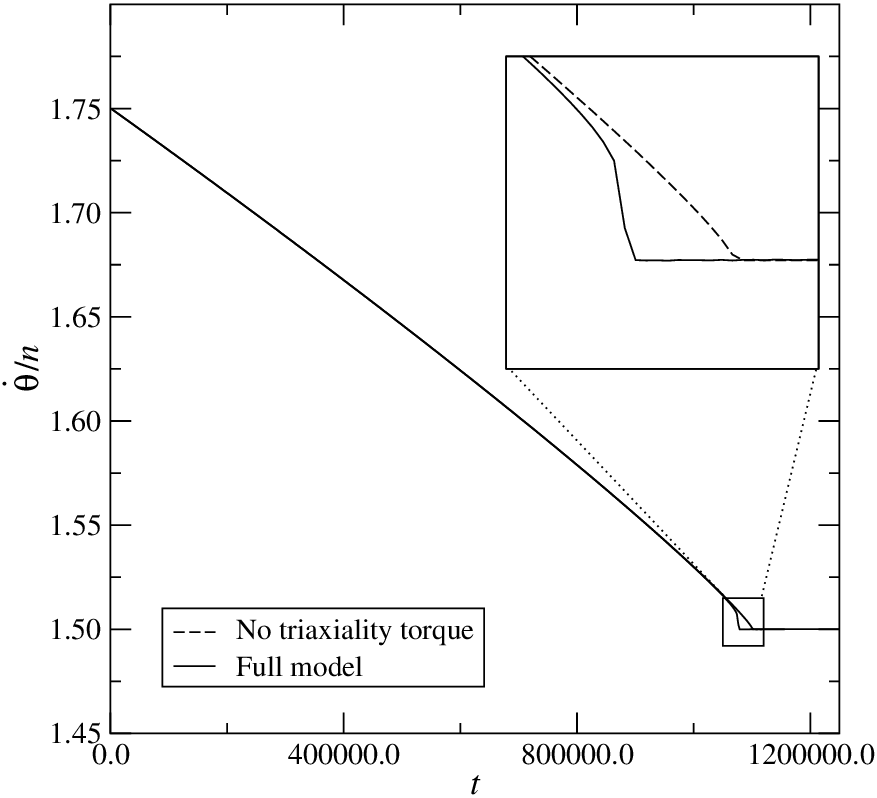}
\caption{Capture can still take place when the triaxial torque is neglected
but the unapproximated tidal torque is used. Here,
$\dot \theta(0) = 1.75 n$ and capture takes place in a time $1.08\times 10^6$
(full model) and $1.10\times 10^6$ (model without triaxial torque).}
\label{figure:6}
\end{figure}

We concentrate first on \eqref{eq:4.9} as an estimate of the
dynamics leading up to capture; we can also use this to estimate the time to
capture in an orbit of a given spin rate. Despite the fact that the triaxial acceleration
has been neglected and the tidal acceleration has been replaced by a 
simple linear approximation, this simplification gives surprisingly good results ---
see Figure \ref{figure:6}, which shows a large timescale comparison, for times 
of order $10^7$, for two different values of the initial spin rate, $\dot \theta(0) = 1.95n$ and $3.25 n$.

Time-to-capture estimates can be made as follows:

\begin{enumerate}
\itemsep0em
\item If $\dot \theta(0) = 1.95n$, then we compute that $a = 1.423\times 10^{-5}$ 
and $b = 1.894\times 10^{-7}$. Neglecting the oscillatory term in the last square brackets in \eqref{eq:4.9} and defining
\begin{equation} \nonumber
R := \dot \theta(0)  -\frac{a}{b} = -24.26 ,
\end{equation}
the time taken for $\dot \theta(t)$ to decay from $\dot \theta(0)$ to $\dot \theta=1.5n$
is estimated to be $b^{-1}\ln[(R-0.45n)/R]$, which leads to the value $2.08\times 10^6$,
to be compared with the numerical value (see Figure \ref{figure:6}) of $2.14\times 10^6$.
\item If $\dot \theta(0) = 3.25n$, then $a = 6.733\times 10^{-6}$,
$b = 2.056\times 10^{-8}$ and $R = -242.70$. In this case, the
estimated time to reach $\dot \theta = 1.5n$ is $8.38\times 10^6$, to be compared with
the numerical result, which is also  $8.38\times 10^6$.
\end{enumerate}

For behaviour on a small timescale, of order 1, again see Figure~\ref{figure:6},
in which, on the right, we plot $\dot \theta(t)$ over a time interval of width 2. Only the order of magnitude of
the approximation is correct, but this is not surprising given that the
approximation completely neglects the triaxial acceleration.

It is interesting to note that the triaxial torque appears to make very little difference during deceleration:
Figure \ref{figure:6} gives evidence that it only has an important role to play very close to capture; see also \cite{FM1,FM2}.

\section{Probability of capture revisited} \label{sec:8}
\setcounter{equation}{0}

Since the seminal paper by Goldreich and Peale \cite{GP}, the probability of capture of a satellite in a resonance $p$:$q$
has usually been studied as the probability $P(p/q)$ for the satellite to be trapped in that resonance
when its rotation velocity $\dot\theta$ approaches the value $pn/q$. However, it may happen that the satellite
never comes close to a given resonance, because it has previously been captured in another one.
Therefore, it may be useful to redefine the probability of capture of a given resonance as
the fraction of initial conditions whose trajectories are attracted by that resonance. Of course the original rotational state 
of the satellite is not even approximately known, so that one has to fix a suitable region in phase space
from which the initial conditions are taken. For the spin-orbit model defined in Section \ref{sec:2},
all the attractors turn out to be contained inside the region $[0,2\pi]\times [-1.5n,4.5n]$, so that, by assuming
an initially prograde satellite (and exploiting the $\pi$-periodicity of the equations of motion by considering
$\theta$ only in the interval $[0,\pi]$), it is reasonable to confine the choice of initial conditions to the region
$\QQ:=[0,\pi]\times [0,4.5n]$. This is essentially what has been done in \cite{BDG2},
where the probability of capture in the $3$:$2$ resonance for Mercury has been estimated to be about $42\%$.
Indeed, even if the initial velocity $\dot\theta$ is much higher, eventually it decreases almost linearly, up to small oscillations
--- as suggested by the analysis performed in Section \ref{sec:4} --- until it enters the region $\QQ$.

Nevertheless, there is no reason why the basins of attraction should be uniformly distributed in phase space.
In fact, Makarov's results \cite{M}, providing an estimated value $P(3/2)=1$, hint that all solutions starting above the $3$:$2$
resonance never reach the attractors contained in the region $\{(\theta,\dot\theta) : \dot\theta<1.5n\}$:
apparently there is a barrier.
Hence, it may be worthwhile to study in more detail the distribution of the basins of attraction. Since they look like
sets of points with no apparent structure, rather than a picture of the basins it is more illuminating to divide the region $\QQ$
into nine strips $\QQ_0,\QQ_1,\ldots,\QQ_8$, each of width $0.5n$ in the velocity direction, and compute the fraction
of the basins of attraction which falls inside each strip.
The results are given in Table \ref{table:5}, where only the prograde resonances have
been considered; we never observed a trajectory with positive initial velocity being attracted by a retrograde resonance.

Table \ref{table:5} confirms the existence of a barrier associated with the $3$:$2$ resonance: 
the trajectories starting with $\dot\theta>1.5n$ cannot reach the attractors below the
$3$:$2$ resonance and, \textit{vice versa}, the trajectories starting with $\dot\theta<1.5n$ cannot reach the attractors above.
Surprisingly, an analogous barrier exists associated with the resonance $1$:$1$.
More precisely, the basins of attraction of the resonances above $3$:$2$  are contained in the region
$\QQ_a:=\{(\theta,\dot\theta) : \dot\theta>1.5n\}$ and, similarly, the basins of attraction of the resonances below
$1$:$1$  are contained in the region $\QQ_b:=\{(\theta,\dot\theta) : \dot\theta<n\}$.
Moreover, all trajectories with initial velocity $\dot\theta\in(n,1.5n)$ are attracted either by the $3$:$2$ resonance
or by the $1$:$1$ resonance. Obviously, we cannot claim that the barriers completely obstruct the passage
of trajectories. In principle, it is possible that holes appear in the barriers;
in any case, crossings of the barriers --- if they occur at all --- appear to be very rare events. 

\begin{table}[!ht]
\centering
\begin{center} 
\setlength\tabcolsep{5pt}
\vskip.5truecm
\vrule
\begin{tabular}{|lr|l|l|l|l|l|l|l|l|l} \hline
$\raisebox{2.2ex}{}\dot\theta(0)/n\;\;\;\in$&&
	$(0,\,0.5]$ & $(0.5,\,1]$ & $(1,\,1.5]$ & $(1.5,\,2]$ & $(2,\,2.5]$ &
	$(2.5,\,3]$ & $(3,\,3.5]$ & $(3.5,\,4]$ & $(4,\,4.5]$\\ \hline
& $m$    & 305  &   36    &  0   &   0   &   0    &  0   &   0   &   0   &   0\\
1:2& $p$ & 0.055&  0.007 &  0   &   0   &   0    &  0   &   0   &   0   &   0\\
& $c$    & 0.006&  0.002 &  0   &   0   &   0    &  0   &   0   &   0   &   0\\ \hline
& $m$    & 5251 &   5493  &  5180  &  0   &    0   &    0   &    0   &    0   &    0\\
1:1& $p$ & 0.945&   0.993 &  0.914 &  0   &    0   &    0   &    0   &    0   &    0\\
& $c$    & 0.006&   0.002 &  0.007 &  0   &    0   &    0   &    0   &    0   &    0\\ \hline
& $m$    & 0   &   0   &   488  &  5432 &  3177 &  2803 &  2939 &  2758 &  2882\\
3:2& $p$ & 0   &   0   &   0.086&  0.990&  0.569&  0.508&  0.519&  0.507&  0.519\\
& $c$    & 0   &   0   &   0.007&  0.003&  0.013&  0.013&  0.013&  0.013&  0.013\\ \hline
& $m$    & 0   &   0   &   0    &  56   &  2361 &  2065 &  2067 &  2013 &  2005\\
2:1& $p$ & 0   &   0   &   0    &  0.010&  0.423&  0.374&  0.365&  0.370&  0.361\\
& $c$    & 0   &   0   &   0    &  0.003&  0.013&  0.013&  0.013&  0.013&  0.013\\ \hline
& $m$    & 0   &   0   &   0   &   0   &   47  &   621  &  523  &  527  &  543\\
5:2& $p$ & 0   &   0   &   0   &   0   &   0.008&  0.113&  0.092&  0.097&  0.098\\
& $c$    & 0   &   0   &   0   &   0   &   0.002&  0.008&  0.008&  0.008&  0.008\\ \hline
& $m$    & 0   &   0   &   0   &   0   &   0   &   28  &   110  &  102  &  93\\
3:1& $p$ & 0   &   0   &   0   &   0   &   0   &   0.005&  0.019&  0.019&  0.017\\
& $c$    & 0   &   0   &   0   &   0   &   0   &   0.002&  0.004&  0.004&  0.003\\ \hline
& $m$    & 0   &   0   &   0   &   0   &   0   &   0   &   23  &   31   &  19\\
7:2& $p$ & 0   &   0   &   0   &   0   &   0   &   0   &   0.004&  0.006&  0.003\\
& $c$    & 0   &   0   &   0   &   0   &   0   &   0   &   0.002&  0.002&  0.002\\ \hline
& $m$    & 0   &   0   &   0   &   0   &   0   &   0   &   0   &   10   &  12\\
4:1& $p$ & 0   &   0   &   0   &   0   &   0   &   0   &   0   &   0.002&  0.002\\
& $c$    & 0   &   0   &   0   &   0   &   0   &   0   &   0   &   0.001&  0.001\\ \hline
Totals & &   5556  & 5529   &5668  & 5488  & 5585   & 5517 &  5662 &  5441 &  5554\\ \hline
\end{tabular}
\hspace{-0.07cm}\vrule
\caption{Estimates of 
the probability that a
solution starting in $\QQ_i := [0,\pi]\times(in/2, (i+1)n/2]$, $i = 0, \ldots, 8$,
ends up in resonance $j$:$2$, this being the resonance with $\dot\theta \approx jn/2$.
For each column, we show $m$, the number of initial conditions that make
the transition; $p$, an estimate of the probability of this happening; and $c$, the 95\%
confidence interval for this probability. The total number of initial
conditions considered for each $\QQ_i$ is shown at the foot of each column,
and the overall total was 50~000.}
\label{table:5}
\end{center}
\end{table}
 
Another interesting point seen in the results of Table \ref{table:5} is that the basin of attraction
of the $3$:$2$ resonance dominates the strip $\QQ_3:=\{(\theta,\dot\theta) : 1.5 n < \dot\theta \le 2n\}$,
the area of its intersection with the strip being about $99\%$ of the whole strip, and becomes more or less uniformly
distributed above the resonance $2$:$1$. Indeed about $51\%$ of the initial conditions in each strip with $\dot\theta>2.5 n$
end up being captured in the $3$:$2$ resonance. Therefore, if we assume the initial velocity $\dot\theta$
to be high enough (as explained above, in practice this means it is sufficient to fix the initial condition above the highest
($4$:$1$) resonance; usually in the literature one takes a value $\dot\theta \approx 4.4 n$;
see for instance \cite{CL,NFME}), then the probability of capture in the $3$:$2$ resonance is more
than $50\%$ --- a value higher than that given in \cite{NFME}, and comparable with that found by
Correia and Laskar \cite{CL} for the CTL model with varying eccentricity.
 
\section{Conclusions} \label{sec:7}
\setcounter{equation}{0}

We have studied the attractors of the spin-orbit model with the realistic tidal torque used in~\cite{NFME}.
First, we have investigated numerically the dynamics of the system: besides the periodic attractors for which the
frequency is locked in a resonance with Mercury's mean motion, quasi-periodic attractors are also detected.
Which attractors arise actually depends on the values of the parameters; in particular quasi-periodic attractors
bifurcate from periodic solutions when the latter become unstable.
For the physical values of the parameters, the main attractor of the Mercury-Sun system is a quasi-periodic
attractor close to the resonance 3:2. This means that, according to the NMFE model,
the librations of the spin rate are quasi-periodic in time in the case of Mercury.

Thereafter, by using a suitable iteration scheme based on perturbation theory for non-smooth systems,
we have provided an explicit analytical expression for the attracting solutions: such expressions,
despite being obtained after a few steps of the iteration and hence being only approximate,
match closely the numerical solutions. So we deduce \textit{a posteriori} that the perturbative approach provides
a reliable description of the dynamics.

However, there are a few aspects which the analysis we have performed does not account for:
\begin{enumerate}
\itemsep0em
\item The quasi-periodic attractor close to the 3:2 resonance has two frequencies: the fast one is 
the mean motion $n$, while the slow one depends on the parameters. We expect the latter to be
slightly different from the frequency $\omega$ obtained by taking into
account just a single harmonic from the Fourier expansion of the triaxial torque --- see \eqref{eq:5.2}.
The second approximation computed in Section \ref{sec:5}
is correct in describing such a phenomenon, but, from a quantitative point of view,
does not provide the right value --- that is to say, the value found numerically (see Section \ref{sec:5.4}).
\item When constructing the second approximation, we have to impose the
condition that a certain integral
vanishes --- see \eqref{eq:5.25a} --- in order to fix the amplitude $C_1$ of the leading term.
This leaves two values: one of them corresponds to the solution which correctly describes
the quasi-periodic attractor found numerically, while the other one has been discarded (see Remark \ref{rmk:5}
in Section \ref{sec:5.3}).
\end{enumerate}

As far as the slow frequency is concerned, the first attempt would be to study the third approximation
and check how the value of the slow frequency changes. However, this is non-trivial
because of the form of the tidal torque, which makes the analytical computations rather tricky:
going beyond the second approximation requires handling equations which involve functions
expressed as integrals, over very long times, of non-smooth functions depending quasi-periodically on time.

As to the discarded second approximation,
we conjecture that it does not correspond to any real solution to the equation because either it disappears
at some higher step of the iteration or the iteration scheme does not converge in such a case --- both possibilities
are likely to be difficult to check in practice.

Another issue that deserves further investigation is the presence of barriers in phase space discussed in Section \ref{sec:8}.
Indeed, it is not obvious why the trajectories cannot cross the resonance $3$:$2$ from above and
the resonance $1$:$1$ from below. Even if the $3$:$2$ resonance corresponds to a quasi-periodic solution
(see Sections \ref{sec:5} and \ref{sec:6}), such a solution is not a KAM torus, which may create
an obstruction in phase space
(in contrast to what happens in the case of the quasi-synchronous resonance arising  in the CTL model);
quite the reverse, it is very localised. Moreover nothing similar happens for the other quasi-periodic attractors,
such as the $2$:$1$ resonance (the second dominant attractor); 
in addition the $1$:$1$ resonance corresponds to a periodic solution.
All this suggests that the occurrence of the barriers is not due to quasi-periodicity. Rather,
it is likely that the phenomenon
is related to the amplitude of the peaks appearing in the tidal torque: indeed the largest peaks correspond
to the $3$:$2$ and $1$:$1$ resonances.


\end{document}